\documentclass[final,1p,times]{elsarticle}
\usepackage{amssymb}
\usepackage{amsmath}
\usepackage{amsthm}
\usepackage[utf8]{inputenc}
\usepackage{marginnote}
\usepackage{amsfonts}
\usepackage{textcomp}
\usepackage{color,soul}
\usepackage{placeins}
\usepackage{tikz}
\usepackage{makecell}
\usepackage{amssymb}
\usepackage{graphicx}
\usepackage{hyperref}
\usepackage[none]{hyphenat}
\usepackage{lipsum}
\newtheorem{theorem}{Theorem}[section]
\newtheorem{proposition}[theorem]{Proposition}
\newtheorem{example}[theorem]{Example}
\newtheorem{lemma}[theorem]{Lemma}
\newtheorem{corollary}[theorem]{Corollary}

\theoremstyle{definition}
\newtheorem{definition}[theorem]{Definition}
\newtheorem{remark}[theorem]{Remark}

\newcommand{\Hn}{\mathcal H_n}
\newcommand{\R}{\mathbb R}
\newcommand{\Ran}{\operatorname{Ran}}
\newcommand{\Span}{\operatorname{span}}

\allowdisplaybreaks
\journal{}

\begin{document}

\begin{frontmatter}

\title{The geometry and classification of positive frames}

\author[ucf]{Shankhadeep Mondal}
\ead{shankhadeep.mondal@ucf.edu}
\author[ucf]{Deguang Han}
\ead{deguang.han@ucf.edu}
\author[tamu]{David Larson}
\ead{d-larson@tamu.edu}
\author[buen]{Ursula M. Molter}
\ead{umolter@dm.uba.ar}

\address[ucf]{Department of Mathematics,
University of Central Florida,
Orlando, FL 32816, USA}

\address[tamu]{Department of Mathematics,
Texas A\&M University,
College Station, TX 77843, USA}
\address[buen]{Department of Mathematics,  University of Buenos Aires, Argentina
}

		\cortext[shan]{Corresponding Author: Shankahdeep Mondal}

\begin{abstract}
We investigate the geometric and operator-theoretic
structure of positive frames and their duals in
finite-dimensional real Hilbert spaces, with particular
emphasis on positive linear dependence, redundancy,
and convex geometry. Building on the classical theory
of positive spanning sets, we characterize equivalence
classes of positive frames through the kernels of their
synthesis operators and obtain parametrizations in
terms of the relative interior of a standard simplex
and suitable subsets of Grassmannians. We establish
that every redundant frame admits a positive dual,
regardless of whether the original frame is positive,
and provide an explicit construction through a
rank-one perturbation of the canonical dual.
Furthermore, we investigate the affine and topological
structure of positive dual frames, proving their
relative openness in the space of all duals and the
absence of extreme points. We characterize dual
frames satisfying prescribed positive relations and
derive criteria for the existence of positive-basis
duals. We also examine semi-positive frames and
their boundary behavior, together with connections
between positive spanning, balancedness, and finite
group representations. Finally, we study the convex
geometry of nonnegative relation spaces and establish
a directional characterization of positive-spanning
robustness under arbitrary vector erasures. These
results provide a unified geometric framework for
understanding positive frames, their dual systems,
and the role of redundancy in constrained
reconstruction.
\end{abstract}


\end{frontmatter}

\section{Introduction}

Redundancy is one of the fundamental features that
distinguishes frames from bases. While a basis provides
a unique linear representation of each vector, a frame
allows multiple representations and, consequently,
considerable freedom in the choice of reconstruction
operators. This flexibility becomes particularly
significant when additional geometric constraints
are imposed on the frame vectors or their duals.
Among such constraints, positive spanning occupies
a distinguished position: it requires that every
vector in the underlying space be generated using
only nonnegative coefficients. The purpose of this
paper is to investigate how this classical geometric
property interacts with the algebraic and operator-
theoretic structure of finite frames, with particular
emphasis on equivalence, duality, redundancy, and
robustness.

The theory of positive spanning sets has a long
history, originating in the study of positive linear
dependence. In his foundational work, Davis
\cite{Davis1954} developed fundamental properties
of positive spanning families and positive bases,
including geometric characterizations and
cardinality restrictions for minimal positive
spanning systems. These ideas subsequently
appeared in the geometric study of convex
polytopes and Gale diagrams \cite{Marcus1984}.
Positive spanning sets have also become important
in derivative-free optimization, where they provide
collections of search directions capable of detecting
descent without explicit derivative information;
see \cite{Conn2009,Audet2017}. 
Dynamical sampling provides another important connection to frame theory through operator-generated systems; see \cite{AldroubiCabrelliKrishtalMolter2026} for a comprehensive survey and \cite{Aldroubi2026} for recent constructions of positive frames with optimal redundancy under erasures.
Related geometric
configurations arise in redundant control systems,
where suitable collections of control directions
are needed to generate prescribed forces and
motions; see, for example,
\cite{Crawford1969,Crawford1968,SanchezPena2000,
Giribet2016}. These developments demonstrate
the usefulness of positive spanning in several
different settings. However, their principal
emphasis is on positive dependence, minimality,
geometric configurations, or applications, rather
than on the structure of dual frames and the
operator-theoretic consequences of positivity.

In finite-dimensional frame theory, every spanning
family is a frame, but not every frame is positively
spanning. Moreover, positivity introduces a
requirement that differs fundamentally from
ordinary linear spanning. A family
$F=\{f_i\}_{i=1}^{N}$ in an $n$-dimensional real
Hilbert space $\mathcal H_n$ is called a
\emph{positive frame} if
\[
\operatorname{span}_{+}(F)
=
\left\{
\sum_{i=1}^{N}c_i f_i:c_i\geq 0
\right\}
=
\mathcal H_n.
\]
The nonnegativity in this definition concerns
the coefficients in a representation, not the
coordinates of the frame vectors themselves.
It is also important to distinguish positive
spanning from reconstruction by a fixed linear
coefficient map. Although a positive frame
allows each vector to be expressed using
nonnegative coefficients, those coefficients
may depend nonlinearly on the vector. In
contrast, a dual frame supplies a prescribed
linear reconstruction formula. Understanding
the relationship between these two forms of
representation is one of the motivations
for the present work.

The classical geometric theory provides a
natural starting point. For a finite spanning
family, positive spanning is equivalent to
the existence of a strictly positive linear
dependence among its vectors, or, equivalently,
to the origin lying in the interior of their
convex hull. In the language of frame
operators, this condition becomes
\[
F\text{ is positive}
\quad\Longleftrightarrow\quad
\ker T_F^*\cap\mathbb R_{++}^{N}
\neq\varnothing,
\]
where $T_F^*$ denotes the synthesis operator.
Although this equivalence follows from
classical positive-spanning theory, its
operator formulation suggests a broader
perspective: the geometry of a positive
frame can be investigated through its
space of linear relations. In particular,
it raises the question of how much of
the structure of a positive frame is
determined by its synthesis kernel.

Our first objective is therefore to study
positive frames under invertible linear
equivalence. Since two ordered frames
are equivalent precisely when their
synthesis operators have the same kernel,
the classification problem can be
reformulated in terms of subspaces of
$\mathbb R^N$ that intersect the strictly
positive orthant. We establish a
parametrization of equivalence classes
of positive frames through these relation
spaces. For frames with excess one, the
parameter space reduces to the relative
interior of a standard simplex, with each
point representing a unique normalized
strictly positive dependence. For
arbitrary excess, we obtain a quotient
description in terms of linearly
independent strictly positive relations,
together with an equivalent Grassmannian
formulation. We also examine the role
of orthogonal projections in constructing
positive frames. These results provide
a geometric classification that depends
on the relation space rather than on
the frame operator alone.

A more substantial question concerns
duality. Given a frame $F$, one may ask
whether its reconstruction freedom can
be used to select a dual frame that
positively spans the ambient space.
This is not the same as asking whether
$F$ itself is positive: the original
frame and its duals can have very
different positive-spanning properties.
We prove that a finite frame admits
a positive dual if and only if it is
redundant. Thus the existence of a
positive dual is governed entirely
by the excess of the original frame.
Our proof is constructive. Starting
from the canonical dual, we use a
strictly positive vector outside the
analysis range to produce a rank-one
perturbation satisfying both the
reconstruction identity and a prescribed
strictly positive linear relation.
This gives an explicit mechanism for
introducing positivity into a dual
system, even when the original frame
is not positively spanning.

Beyond existence, we investigate the
geometry of the family
$\mathcal D_{+}(F)$ of all positive
duals of a fixed frame. We establish
its relative openness within the
affine space of dual frames and
examine the perturbations that
preserve positivity and reconstruction
simultaneously. In particular, we
show that $\mathcal D_{+}(F)$ has
no extreme points: every positive
dual can be expressed as the midpoint
of two distinct positive duals.
We also study dual frames satisfying
prescribed positive relations and
the more restrictive problem of
obtaining a dual that is a positive
basis. The latter leads to conditions
involving complements of the analysis
range and prescribed coordinate sign
patterns. Together, these results
describe how positive spanning
restricts, but does not eliminate,
the freedom available in dual-frame
reconstruction.

The geometry becomes more delicate
when strict positivity is relaxed.
Motivated by this distinction, we
consider semi-positive frames,
which admit nontrivial nonnegative
linear relations with at least one
vanishing coefficient. We investigate
their approximation by positive
frames and the corresponding
boundary behavior of positive duals.
This provides a connection between
strictly positive relations and
degenerate nonnegative dependencies,
and clarifies how positivity may
be lost under limiting processes.
The resulting picture is useful
for understanding the topology
of families of positive frames
and their dual systems.

We further examine positive frames
arising from finite group
representations. Group-generated
families possess additional
symmetry, and their linear relations
reflect the structure of the
underlying representation. We
investigate the relationship between
balancedness, linear connectedness,
and positive spanning, using the
averaging operator associated with
a finite group. In particular, the
presence or absence of invariant
vectors provides information about
balanced orbit families and their
positivity. This places the
geometric theory of positive
frames in a structured setting
where algebraic symmetry can
be used to identify positive
dependencies.

Finally, we study positive
dependence from the viewpoint
of convex geometry and
robustness. Minimal positive
dependencies and the extreme
points of the polytope of
normalized nonnegative
relations provide information
about the internal structure
of positive frames. We also
consider the effect of deleting
frame vectors. Unlike ordinary
spanning, positive spanning
requires every nonzero
direction to have a strictly
positive inner product with
at least one frame vector.
This leads to the directional
quantity
\[
d_{+}(F)
=
\min_{\|x\|=1}
\#\bigl\{
i:\langle x,f_i\rangle>0
\bigr\}.
\]
We show that a positive frame
remains positively spanning
after the deletion of any
$s$ vectors if and only if
$d_{+}(F)\geq s+1$. We
also investigate related
questions of exactness,
cardinality, and approximation
by positive Parseval frames.
These results connect the
combinatorial structure of
positive dependencies with
the stability of positive
spanning under perturbations
and erasures.

The paper is organized as follows.
Section~2 recalls the necessary background on finite
frames, positive spanning sets, and convex geometry.
Section~3 develops geometric characterizations of
positive frames, their equivalence classes through
synthesis kernels, simplex and Grassmannian
parametrizations, and the structure of positive bases.
Section~4 investigates positive dual frames, including
their existence, construction, and geometric properties.
Section~5 examines the convex geometry of positive
duals, their boundary behavior, and the approximation
of semi-positive frames.
Section~6 studies positive frames generated by finite
group representations, with particular emphasis on
balancedness and linear connectedness.
Section~7 explores the convex geometry of nonnegative
relations, positive-spanning robustness under erasures,
exactness, cardinality bounds, and positive Parseval
approximation. Finally, Section~8 presents concluding remarks and
directions for future research.

\section{Preliminaries}\label{sec:preliminaries}

Throughout this paper, $\mathcal H_n$ denotes an
$n$-dimensional real Hilbert space, where $n\geq 1$.
All vector spaces and operators are considered over
$\mathbb R$. We work with finite indexed families
$F=\{f_i\}_{i=1}^N\subset\mathcal H_n$, where
repetitions and zero vectors are permitted unless
otherwise specified. We recall the basic notions
of finite frame theory, positive spanning, and
convex geometry that will be used throughout
the paper.

\subsection{Finite frames and duality}

A family $F=\{f_i\}_{i=1}^N$ is called a
\emph{frame} for $\mathcal H_n$ if there exist
constants $0<A\leq B<\infty$ such that
\[
A\|x\|^2\leq
\sum_{i=1}^N|\langle x,f_i\rangle|^2
\leq B\|x\|^2,
\qquad x\in\mathcal H_n.
\]
In finite dimensions, this is equivalent to
$\operatorname{span}(F)=\mathcal H_n$.
The analysis operator $T_F:\mathcal H_n\to
\mathbb R^N$ and the synthesis operator
$T_F^*:\mathbb R^N\to\mathcal H_n$ are
defined by
\[
T_Fx=\sum_{i=1}^N\langle x,f_i\rangle e_i,
\qquad
T_F^*c=\sum_{i=1}^Nc_if_i,
\]
where $\{e_i\}_{i=1}^N$ is the standard
orthonormal basis of $\mathbb R^N$ and
$c=(c_i)_{i=1}^N$. The associated frame
operator is $S_F=T_F^*T_F$, which is
positive, self-adjoint, and invertible.

A frame is called \emph{tight} if
$S_F=AI_{\mathcal H_n}$ for some $A>0$,
and \emph{Parseval} if $S_F=I_{\mathcal H_n}$.
The canonical dual of $F$ is
$\widetilde F=\{S_F^{-1}f_i\}_{i=1}^N$.

A frame $G=\{g_i\}_{i=1}^N$ is called
a \emph{dual frame} of $F$ if
\[
x=\sum_{i=1}^N\langle x,f_i\rangle g_i,
\qquad x\in\mathcal H_n,
\]
or equivalently, $T_G^*T_F=I_{\mathcal H_n}$.
We denote by $\mathcal D(F)$ the collection
of all dual frames of $F$.

Recall that every dual frame of $F$ has
the representation
\[
g_i=S_F^{-1}f_i+h_i,
\qquad 1\leq i\leq N,
\]
where $H=\{h_i\}_{i=1}^N$ satisfies
$T_H^*T_F=0$. Consequently,
\[
\mathcal D(F)=
\left\{
\widetilde F+H:T_H^*T_F=0
\right\}.
\]
Thus $\mathcal D(F)$ is an affine
subspace of $\mathcal H_n^N$. The \emph{relation space} of $F$ is
$K_F=\ker T_F^*$. Since $T_F^*$ is
surjective, the rank-nullity theorem gives
\[
\dim K_F=N-n,
\qquad
K_F=(\operatorname{Ran}T_F)^\perp.
\]
The integer $N-n$ is called the
\emph{excess} of $F$. In particular,
$F$ is redundant if and only if
$K_F\neq\{0\}$.

We shall frequently use the orthogonal
projection onto $\operatorname{Ran}T_F$,
given by
\[
P_F=T_FS_F^{-1}T_F^*.
\]
Consequently, $I_{\mathbb R^N}-P_F$
is the orthogonal projection onto
$K_F$. These identities provide
a convenient connection between
the reconstruction operators and
the linear relations among frame
vectors.

\subsection{Positive spanning sets and positive bases}

For a finite family
$F=\{f_i\}_{i=1}^N\subset\mathcal H_n$,
its \emph{positive span} is defined by
\[
\operatorname{span}_+(F)
=
\left\{
\sum_{i=1}^Nc_if_i:c_i\geq0
\right\}.
\]
A frame $F$ is called a
\emph{positive frame} if
$\operatorname{span}_+(F)=\mathcal H_n$.
A dual frame $G$ of $F$ is called a
\emph{positive dual frame} if $G$
is itself a positive frame. We write
\[
\mathcal D_+(F)=
\{G\in\mathcal D(F):
\operatorname{span}_+(G)=\mathcal H_n\}.
\]

A family $F$ is said to admit a
\emph{nonnegative relation} if
there exist coefficients $c_i\geq0$,
not all zero, such that
$\sum_{i=1}^Nc_if_i=0$.
If every coefficient is strictly
positive, the relation is called
a \emph{strictly positive relation}.
For convenience, we write
\[
\mathbb R_{\geq0}^N
=\{c\in\mathbb R^N:c_i\geq0\},
\qquad
\mathbb R_{++}^N
=\{c\in\mathbb R^N:c_i>0\}.
\]
Thus nonnegative and strictly positive
relations correspond, respectively,
to nonzero vectors in
$K_F\cap\mathbb R_{\geq0}^N$
and vectors in $K_F\cap\mathbb R_{++}^N$.

The theory of positive linear dependence
was developed systematically by Davis
\cite{Davis1954}. We recall a classical
result that will be used in our analysis. A positive frame $F$ is called a
\emph{positive basis} if no proper
subfamily of $F$ positively spans
$\mathcal H_n$. Thus a positive basis
is a minimal positive spanning family,
which need not be a basis in the usual
linear-algebraic sense. The following
cardinality restriction is another
classical result of Davis.

\begin{theorem}[Davis \cite{Davis1954}]
\label{thm:davis-cardinality}
Let $F=\{f_i\}_{i=1}^N$ be a positive
basis for $\mathcal H_n$. Then
\[
n+1\leq N\leq2n.
\]
\end{theorem}

Positive bases consisting of $n+1$
and $2n$ vectors are called
\emph{minimal positive bases} and
\emph{maximal positive bases},
respectively. These configurations
also arise in derivative-free
optimization, where positive spanning
sets are used to construct suitable
search directions; see
\cite{Conn2009,Audet2017,KoldaLewisTorczon2003}.

The following standard examples
illustrate the distinction between
positive frames and positive bases.

\begin{example}\label{ex:positive-bases}
Let $\{e_i\}_{i=1}^n$ be the standard
orthonormal basis of $\mathbb R^n$.

\begin{enumerate}
\item[(i)]
The family
\[
F_1=
\left\{
e_1,\ldots,e_n,
-\sum_{i=1}^ne_i
\right\}
\]
is a minimal positive basis
consisting of $n+1$ vectors.

\item[(ii)]
The family
\[
F_2=
\{e_1,-e_1,\ldots,e_n,-e_n\}
\]
is a maximal positive basis
consisting of $2n$ vectors.

\item[(iii)]
The family
\[
F_3=
\{(1,0),(0,1),(-1,-1),(1,1)\}
\]
is a positive frame for $\mathbb R^2$,
but not a positive basis, since its
first three vectors already
positively span $\mathbb R^2$.
\end{enumerate}
\end{example}

\subsection{Convex geometry and positive relations}

For a finite family $F=\{f_i\}_{i=1}^N$,
its convex hull is
\[
\operatorname{conv}(F)
=
\left\{
\sum_{i=1}^N\alpha_if_i:
\alpha_i\geq0,\ 
\sum_{i=1}^N\alpha_i=1
\right\}.
\]
We use $\operatorname{int}(C)$,
$\operatorname{relint}(C)$, and
$\operatorname{Ext}(C)$ to denote
the interior, relative interior,
and set of extreme points of
a convex set $C$, respectively.

For a cone $C\subseteq\mathcal H_n$,
its negative polar cone is defined by
\[
C^\circ=
\{x\in\mathcal H_n:
\langle x,y\rangle\leq0
\text{ for all }y\in C\}.
\]
In particular, for the finitely
generated cone
$C=\operatorname{span}_+(F)$,
we have
\[
C^\circ=
\{x\in\mathcal H_n:
\langle x,f_i\rangle\leq0
\text{ for all }i\}.
\]
Since finitely generated cones are
closed and convex, the bipolar
theorem yields
\[
\operatorname{span}_+(F)=\mathcal H_n
\quad\Longleftrightarrow\quad
C^\circ=\{0\}.
\]
This gives a convex-analytic
interpretation of \cite{Davis1954}.

Another useful classical fact is
the supporting hyperplane theorem:
if a convex set has nonempty
interior and the origin does not
belong to that interior, then
there exists a nonzero linear
functional separating the origin
from the interior. In the present
finite-dimensional setting, such
functionals can be represented
by inner products.

We shall also use the following
standard consequence of the
separation theorem.

\begin{lemma}\label{lem:convex-hull-separation}
Let $F=\{f_i\}_{i=1}^N$ be a
spanning family for $\mathcal H_n$.
Then
$0\in\operatorname{int}\operatorname{conv}(F)$
if and only if, for every nonzero
$x\in\mathcal H_n$, there exist
indices $i,j$ such that
\[
\langle x,f_i\rangle>0,
\qquad
\langle x,f_j\rangle<0.
\]
\end{lemma}

\begin{proof}
Suppose that
$0\in\operatorname{int}\operatorname{conv}(F)$.
For every nonzero $x\in\mathcal H_n$,
there exists $\varepsilon>0$ such that
$\pm\varepsilon x\in\operatorname{conv}(F)$.
Taking inner products with $x$
shows that some frame vector has
positive inner product with $x$,
while another has negative inner
product with $x$.

Conversely, suppose the stated
sign condition holds. If the
origin were not in the interior
of $\operatorname{conv}(F)$, the
separation theorem would provide
a nonzero $x\in\mathcal H_n$
such that
$\langle x,f_i\rangle\geq0$
for every $i$, contradicting
the assumed condition.
\end{proof}

The preceding lemma is a
convenient form of the classical
positive-spanning criterion.
In Section~3, we shall relate
these geometric properties
to strictly positive vectors
in the synthesis kernel.

For later use, we also introduce
the set of normalized
nonnegative relations
\[
\Omega_F=
\left\{
c\in K_F\cap\mathbb R_{\geq0}^N:
\sum_{i=1}^Nc_i=1
\right\}.
\]
This is a compact convex polytope,
possibly empty. Its elements
encode nonnegative linear
dependencies among the frame
vectors. We shall investigate
its geometric structure and
extreme points in a later section.

Finally, we recall that a
convex subset $C$ of an affine
space has a relative interior
taken with respect to
$\operatorname{aff}(C)$.
The distinction between
ordinary and relative
interiors will be important
when studying families of
dual frames, since
$\mathcal D(F)$ is generally
a proper affine subspace
of $\mathcal H_n^N$.

\section{Geometry and Classification of Positive Frames}

Throughout this section, let $\mathcal{H}_n$ be an
$n$-dimensional real Hilbert space and let
$F=\{f_i\}_{i=1}^N$ be a frame for $\mathcal{H}_n$. 
In this section, we investigate the geometric and
operator-theoretic structure of positive frames,
with particular emphasis on their characterization
and classification under invertible linear
transformations. We begin by relating positive
spanning to strictly positive linear dependencies
and the interior geometry of convex hulls.
Building on these classical characterizations,
we develop a relation-space approach to the
classification of positive frames. In particular,
we show that their equivalence classes can be
completely described through the kernels of
their synthesis operators. This leads to a
parametrization by the relative interior of
a standard simplex in the excess-one case
and, more generally, by suitable subsets of
Grassmannians. We also investigate constructions
arising from orthogonal projections and demonstrate
that positive Parseval frames with identical
frame operators need not be equivalent. These
results provide a geometric description of
positive frames in which the synthesis kernel,
rather than the frame operator alone, plays
the central role.

\subsection{\bf{Geometric Characterizations of Positive Frames}} 

We begin by recalling equivalent geometric
characterizations of positive frames in terms
of strictly positive linear dependencies and
the interior of their convex hulls. We then
interpret these conditions through the kernel
of the synthesis operator and establish
some related properties that will be useful
in the subsequent classification of positive
frames.

The following result provides equivalent
characterizations of positive frames in terms of
strictly positive linear dependence and the
geometry of the convex hull. 

\begin{proposition}\label{thm:positive-characterization}
Let $F=\{f_i\}_{i=1}^{N}$ be a frame for
$\mathcal{H}_n$. Then the following statements
are equivalent:
\begin{enumerate}
\item[(i)] $F$ is a positive frame.
\item[(ii)] There exist scalars
$\lambda_1,\ldots,\lambda_N>0$ such that $\sum_{i=1}^{N}\lambda_i f_i=0.$
\item[(iii)] $0\in\operatorname{int}
\bigl(\operatorname{conv}(F)\bigr)$.
\end{enumerate}
\end{proposition}

\begin{proof}
We prove the implications
$(i)\Rightarrow(ii)\Rightarrow(iii)\Rightarrow(i)$.

$(i)\Rightarrow(ii)$.
Suppose that $F$ is a positive frame.
Since $-f_i\in\operatorname{span}_{+}(F)$
for every $1\leq i\leq N$, there exist
$a_{ij}\geq 0$ such that
\[
-f_i=\sum_{j=1}^{N}a_{ij}f_j.
\]
Summing over $i$ and rearranging yields
\[
\sum_{j=1}^{N}
\left(1+\sum_{i=1}^{N}a_{ij}\right)f_j=0.
\]
Setting
$\lambda_j=1+\sum_{i=1}^{N}a_{ij}>0$
establishes (ii).

$(ii)\Rightarrow(iii)$.
Suppose that
$\sum_{i=1}^{N}\lambda_i f_i=0$
for some $\lambda_i>0$.
Normalizing the coefficients, we obtain
\[
0=\sum_{i=1}^{N}\alpha_i f_i,
\qquad
\alpha_i=\frac{\lambda_i}
{\sum_{j=1}^{N}\lambda_j}>0,
\qquad
\sum_{i=1}^{N}\alpha_i=1.
\]
Thus $0\in\operatorname{conv}(F)$. Assume that $0$ does not belong to the
interior of $\operatorname{conv}(F)$.
By the supporting hyperplane theorem,
there exists a nonzero $x\in\mathcal{H}_n$
such that $\langle x,f_i\rangle\geq 0$
for every $i$. Consequently,
\[
0=
\left\langle x,
\sum_{i=1}^{N}\lambda_i f_i
\right\rangle
=
\sum_{i=1}^{N}\lambda_i
\langle x,f_i\rangle.
\]
Since every $\lambda_i>0$, it follows that
$\langle x,f_i\rangle=0$ for all $i$.
As $F$ spans $\mathcal{H}_n$, we obtain
$x=0$, a contradiction. Hence
$0\in\operatorname{int}(\operatorname{conv}(F))$.

$(iii)\Rightarrow(i)$.
Suppose that
$0\in\operatorname{int}(\operatorname{conv}(F))$.
Then there exists $\varepsilon>0$ such that
$B(0,\varepsilon)\subseteq\operatorname{conv}(F)$.
For any nonzero $x\in\mathcal{H}_n$, we have
\[
\frac{\varepsilon x}{2\|x\|}
\in\operatorname{conv}(F).
\]
Therefore, there exist $c_i\geq 0$ with
$\sum_{i=1}^{N}c_i=1$ such that
\[
\frac{\varepsilon x}{2\|x\|}
=\sum_{i=1}^{N}c_i f_i.
\]
It follows that
\[
x=\sum_{i=1}^{N}
\frac{2\|x\|c_i}{\varepsilon}f_i
\in\operatorname{span}_{+}(F).
\]
Since $0$ also belongs to the positive span,
we conclude that
$\operatorname{span}_{+}(F)=\mathcal{H}_n$.
\end{proof}

The preceding theorem also admits a useful
operator-theoretic formulation.

\begin{remark}
Let $\mathbb{R}_{++}^{N}
=
\{c=(c_i)_{i=1}^{N}\in\mathbb{R}^{N}:
c_i>0 \text{ for all }i\}.$ Then Proposition~\ref{thm:positive-characterization}
implies that
\[
F\text{ is a positive frame}
\quad\Longleftrightarrow\quad
\ker T_F^*\cap\mathbb{R}_{++}^{N}
\neq\varnothing.
\]
Thus the positivity of a frame is completely
determined by the existence of a strictly
positive vector in its relation space
$\ker T_F^*$.
\end{remark}

\begin{example}
Consider the family $F=\{f_i\}_{i=1}^5$ in $\mathbb R^3$, where
\[
\begin{aligned}
f_1&=(1,0,0), \quad
f_2=(0,1,0), \quad
f_3=(0,0,1),\quad
f_4&=(-1,-1,0), \quad
f_5=(-1,0,-1).
\end{aligned}
\]
Since $f_1,f_2,f_3$ form the standard basis of
$\mathbb R^3$, the family $F$ is a frame with excess two.
Its synthesis operator is represented by
\[
T_F^*=
\begin{pmatrix}
1&0&0&-1&-1\\
0&1&0&-1&0\\
0&0&1&0&-1
\end{pmatrix}.
\]
A direct computation yields
\[
\ker T_F^*
=
\{(s+t,s,t,s,t):s,t\in\mathbb R\}.
\]
In particular, $(2,1,1,1,1)\in\ker T_F^*$,
and hence $F$ is a positive frame by
Proposition~\ref{thm:positive-characterization}.

Furthermore, the normalized strictly positive
relations are precisely
\[
\left\{
\left(\frac13,s,\frac13-s,s,\frac13-s\right):
0<s<\frac13
\right\}.
\]
Thus, unlike the excess-one case, the normalized
positive relation is not unique.

Notice also that
\[
f_1+f_2+f_4=0,
\qquad
f_1+f_3+f_5=0.
\]
These relations correspond to two distinct
positively dependent subfamilies, neither of
which positively spans $\mathbb R^3$.
However, their union positively spans
$\mathbb R^3$.

Finally, the subfamily obtained by deleting
$f_2$ still spans $\mathbb R^3$, but it is
not positive. Indeed, every relation supported
on the remaining indices must satisfy $s=0$,
and therefore cannot have all coefficients
strictly positive.
\end{example}

We observe that the convex hull of a positive frame
satisfies a useful scaling property. Indeed, let
$F=\{f_i\}_{i=1}^N$ be a positive frame for
$\mathcal H_n$. Since
$0\in\operatorname{int}\operatorname{conv}(F)$,
there exists $r>0$ such that
$B(0,r)\subseteq\operatorname{conv}(F)$.
Now, for any $x\in\mathcal H_n$, choose
$t>\|x\|/r$. Then
\[
\frac{x}{t}\in B(0,r)
\subseteq\operatorname{conv}(F),
\]
which implies that
$x\in t\operatorname{conv}(F)$.
Thus, every vector in $\mathcal H_n$
belongs to a suitable positive scaling
of the convex hull of $F$.

\begin{proposition}\label{lem:positive-relation-shift}
Let $F=\{f_i\}_{i=1}^N$ be a positive
frame for $\mathcal H_n$, and let
$\{e_i\}_{i=1}^N$ be the standard
orthonormal basis of $\mathbb R^N$.
Suppose that $\sum_{i=1}^N\lambda_i f_i=0,
\quad \lambda_i>0,$
and set $u=\sum_{i=1}^N\lambda_i e_i$.
Then, for every $x\in\mathcal H_n$,
there exists $t\geq0$ such that
\[
T_Fx+tu\in
\operatorname{span}_+\{e_1,\ldots,e_N\}.
\]
\end{proposition}

\begin{proof}
Let $x\in\mathcal H_n$ be arbitrary. We then have
\[
T_Fx+tu
=\sum_{i=1}^N
\bigl(\langle x,f_i\rangle
+t\lambda_i\bigr)e_i.
\]
Since $\lambda_i>0$ for every $i$,
we may choose
\[
t\geq
\max\left\{
0,\max_{1\leq i\leq N}
\frac{-\langle x,f_i\rangle}{\lambda_i}
\right\}.
\]
Then
$\langle x,f_i\rangle+t\lambda_i\geq0$
for every $i$. Consequently, $T_Fx+tu\in
\operatorname{span}_+\{e_1,\ldots,e_N\},$
as required.
\end{proof}

\subsection{\bf{Equivalence Classes of Positive Frames}}\label{sec:equivalence}
In the preceding section, we characterized positive frames
through strictly positive linear relations among their
vectors. We now investigate the structure of positive frames
under invertible linear transformations. Our main objective
is to describe their equivalence classes in terms of the
kernels of their synthesis operators. In particular, we
show that these classes admit a natural parametrization
by certain subspaces of $\mathbb R^N$ intersecting the
strictly positive orthant. The case of frames with excess
one is considered separately, where the corresponding
parameter space has a particularly simple geometric
description.

Two frames $F=\{f_i\}_{i=1}^{N}
\quad\text{and}\quad
G=\{g_i\}_{i=1}^{N}$ for $\Hn$ are said to be \emph{equivalent} if there exists an invertible operator $U:\Hn\longrightarrow\Hn$ such that
$g_i=Uf_i,
\quad i=1,\ldots,N.$
In this case we write $F\sim G.$ For finite frames, equivalence is completely determined by the space of
linear relations. More precisely,
\[
F\sim G
\quad\Longleftrightarrow\quad
\ker T_F^*=\ker T_G^*.
\]
Indeed, $g_i=Uf_i$ with $U$ invertible implies $T_G^*=UT_F^*$,
so the kernels agree. Conversely, if the kernels agree, the rule
$U(T_F^*c)=T_G^*c$ is well defined and is a linear bijection, with
$Uf_i=g_i$. Thus equivalence preserves positivity by
Theorem~\ref{thm:positive-characterization}.

This observation, together with the existence of strictly positive
relations, gives a natural quotient parametrization of positive frames. We first consider positive frames with excess one.
Let $\mathbb R_+^{n+1}
=
\left\{
\lambda=(\lambda_i)_{i=1}^{n+1}\in\mathbb R^{n+1}:
\lambda_i\geq 0
\right\}
\;\text{and}\;
\mathbb R_{++}^{n+1}
=
\left\{
\lambda=(\lambda_i)_{i=1}^{n+1}\in\mathbb R^{n+1}:
\lambda_i>0
\right\}$
denote the nonnegative and strictly positive orthants
of $\mathbb R^{n+1}$, respectively. Define
\[
\omega
=
\operatorname{conv}\{e_1,\ldots,e_{n+1}\}
=
\left\{
\lambda\in\mathbb R_+^{n+1}:
\sum_{i=1}^{n+1}\lambda_i=1
\right\},
\]
where $\{e_i\}_{i=1}^{n+1}$ is the standard basis of
$\mathbb R^{n+1}$. Its relative interior is given by
\[
\mathring{\omega}
=
\left\{
\lambda\in\mathbb R_{++}^{n+1}:
\sum_{i=1}^{n+1}\lambda_i=1
\right\}.
\]
The following proposition characterizes the equivalence
classes of positive frames with excess one in terms of
the relative interior of the standard simplex.

\begin{proposition}\label{prop:positive-frame-simplex}
The equivalence classes of ordered positive frames
consisting of $n+1$ vectors in an $n$-dimensional
real Hilbert space $\mathcal H_n$ are in one-to-one
correspondence with the points of $\mathring{\omega}$.
\end{proposition}

\begin{proof}
Let $F=(f_i)_{i=1}^{n+1}$ be positive. Since $T_F^*$ is surjective,
its kernel has dimension one. By
Theorem~\ref{thm:positive-characterization}, this kernel contains a
strictly positive vector $a=(a_i)$. Thus
\[
 \ker T_F^*=\Span\{a\},\qquad
 \lambda_i=\frac{a_i}{\sum_{j=1}^{n+1}a_j}
\]
defines an element $\lambda\in\mathring\omega$ with
$\sum_i\lambda_if_i=0$. This normalized positive relation is unique:
any other such vector $\mu$ lies in the same one-dimensional kernel,
so $\mu=t\lambda$ for some $t>0$, and summing the coordinates gives
$t=1$.

Define $\Phi([F])=\lambda$. Equivalent frames have the same synthesis
kernel, so this map is well defined. If $\Phi([F])=\Phi([G])$, then
$\ker T_F^*=\Span\{\lambda\}=\ker T_G^*$; the kernel criterion above
implies $F\sim G$. Hence $\Phi$ is injective.

To prove surjectivity, fix $\lambda\in\mathring\omega$, put
$K=\Span\{\lambda\}$, and let $P$ be the orthogonal projection onto
$K^\perp\subset\R^{n+1}$. Define $f_i=Pe_i$. The family $(Pe_i)$
spans the $n$-dimensional real Hilbert space $K^\perp$, and
\[
 \sum_{i=1}^{n+1}\lambda_iPe_i=P\lambda=0.
\]
It is therefore positive. Its synthesis operator is $P$, regarded as a
map into $K^\perp$, with kernel $K$. Consequently its normalized
positive relation is precisely $\lambda$. Finally, choose an isometric
isomorphism from $K^\perp$ onto $\Hn$ to obtain a representative in the
fixed ambient space. This proves surjectivity and the asserted bijection.
\end{proof}

\begin{remark}\label{rem:canonical-positive-representative}
For each $\lambda\in\mathring{\omega}$, one may obtain a concrete
representative of the corresponding equivalence class without using a
quotient space. Consider
\[
\lambda^\perp
=
\left\{
x\in\mathbb R^{n+1}:
\langle x,\lambda\rangle=0
\right\},
\]
which is an $n$-dimensional subspace of $\mathbb R^{n+1}$. Let
$P_{\lambda^\perp}$
denote the orthogonal projection onto $\lambda^\perp$, and define
$f_i^\lambda=P_{\lambda^\perp}e_i, \qquad i=1,\ldots,n+1.$
Then
\[
\sum_{i=1}^{n+1}\lambda_i f_i^\lambda
=
P_{\lambda^\perp}
\left(
\sum_{i=1}^{n+1}\lambda_i e_i
\right)
=
P_{\lambda^\perp}\lambda
=
0.
\]
Thus
\[
\left\{
P_{\lambda^\perp}e_i
\right\}_{i=1}^{n+1}
\]
is a positive frame for $\lambda^\perp$ and provides a canonical
representative of the equivalence class corresponding to
$\lambda\in\mathring{\omega}$.
\end{remark}

\begin{remark}
The preceding proposition concerns \emph{ordered} frames. If
permutations of the frame vectors are also included in the equivalence
relation, then the corresponding parameter space is the quotient
$\mathring{\omega}/S_{n+1},$
where the symmetric group $S_{n+1}$ acts on $\mathring{\omega}$ by
permuting coordinates.
\end{remark}

For the general parametrization, let $k=N-n$ and use the strictly positive orthant $\R_{++}^N$. Define
\[
\mathfrak A_k
=
\left\{
(\lambda^{(1)},\ldots,\lambda^{(k)})\in(\R_{++}^N)^k:
\lambda^{(1)},\ldots,\lambda^{(k)}
\text{ are linearly independent}
\right\}.
\]

On $\mathfrak A_k$, define the relation $\approx$ by
\[
(\lambda^{(1)},\ldots,\lambda^{(k)})\approx(\mu^{(1)},\ldots,\mu^{(k)})
\]
if and only if
$\Span\{\lambda^{(1)},\ldots,\lambda^{(k)}\} = \Span\{\mu^{(1)},\ldots,\mu^{(k)}\}.$

We first consider the construction of positive frames
through orthogonal projections. 
Let $\{e_i\}_{i=1}^N$ be an orthonormal basis
of $\mathbb R^N$, and let
$x=\sum_{i=1}^N\lambda_i e_i$ be a nonzero
vector with $\lambda_i\geq0$ for every $i$.
Denote by $P$ the orthogonal projection
onto $x^\perp$. We now characterize when
the family $\{Pe_i\}_{i=1}^N$ forms a positive
frame for $x^\perp$ in terms of the
coefficients of $x$.

\begin{proposition}\label{prop:projection-positive-frame}
The family $\{Pe_i\}_{i=1}^N$ is a positive
frame for $x^\perp$ if and only if
$\lambda_i>0$ for every $i$.
\end{proposition}

\begin{proof}
Since $P$ is the orthogonal projection
onto $x^\perp$, the family
$\{Pe_i\}_{i=1}^N$ spans $x^\perp$.

Suppose first that $\lambda_i>0$
for every $i$. Since $Px=0$, we have
\[
\sum_{i=1}^N\lambda_i Pe_i
=P\left(\sum_{i=1}^N\lambda_i e_i\right)
=Px=0.
\]
Thus $\{Pe_i\}_{i=1}^N$ admits a
strictly positive linear relation.
Since it spans $x^\perp$, it is a
positive frame.

Conversely, suppose that
$\{Pe_i\}_{i=1}^N$ is a positive frame
for $x^\perp$. Then there exist
$c_i>0$ such that
\[
\sum_{i=1}^Nc_iPe_i=0.
\]
Consequently,
\[
P\left(\sum_{i=1}^Nc_ie_i\right)=0,
\]
which implies
$\sum_{i=1}^Nc_ie_i\in\ker P
=\operatorname{span}\{x\}$.
Hence there exists $t\neq0$ such that
\[
\sum_{i=1}^Nc_ie_i
=t\sum_{i=1}^N\lambda_i e_i.
\]
By the linear independence of
$\{e_i\}_{i=1}^N$, we obtain
$c_i=t\lambda_i$ for every $i$.
Since $c_i>0$ and $\lambda_i\geq0$,
it follows that $\lambda_i>0$
for every $i$.
\end{proof}

\begin{example}
We construct positive Parseval frames for
$\mathbb R^3$ that are not equivalent despite
having the same frame operator. We also
illustrate that identical vector norms and
the Parseval property do not guarantee
positive spanning.

Let
\[
u=\frac{1}{2}(1,1,1,1)^T,
\qquad
v=\frac{1}{2}(1,1,-1,-1)^T,
\]
and denote by $P_u$ and $P_v$ the orthogonal
projections onto $u^\perp$ and $v^\perp$,
respectively. Consider the families
\[
F=\{P_u e_i\}_{i=1}^4,
\qquad
G=\{P_v e_i\}_{i=1}^4.
\]
Both families are Parseval frames for their
respective three-dimensional subspaces.
Moreover,
\[
\|P_u e_i\|^2
=
\|P_v e_i\|^2
=
\frac{3}{4},
\qquad 1\leq i\leq 4.
\]

Since $u$ has strictly positive coordinates,
$F$ is a positive frame. However, $G$ is
not positive because its synthesis kernel
is $\operatorname{span}\{v\}$, which contains
no strictly positive vector.

Now let
\[
w=\frac{1}{\sqrt{10}}(1,1,2,2)^T
\]
and consider the positive Parseval frame
$H=\{P_w e_i\}_{i=1}^4$ for $w^\perp$.
Its synthesis kernel is
$\operatorname{span}\{w\}$, so $H$ is positive.

Although $F$ and $H$ are both positive
Parseval frames with four vectors in
three-dimensional Hilbert spaces, they
are not equivalent. Indeed,
\[
\ker T_F^*=\operatorname{span}\{u\},
\qquad
\ker T_H^*=\operatorname{span}\{w\},
\]
and these subspaces are distinct. Furthermore, their normalized strictly
positive relations are $\lambda_F=
\left(\frac14,\frac14,\frac14,\frac14\right) \quad and \quad \lambda_H=
\left(\frac16,\frac16,\frac13,\frac13\right),$
respectively. Thus Proposition~\ref{prop:positive-frame-simplex}
distinguishes their equivalence classes,
even though both frames are Parseval. This example demonstrates that the
equivalence classification of positive
frames is governed by their relation
spaces rather than by the Parseval
property alone.
\end{example}

\subsection{\bf{Higher Excess and Grassmannian Parametrization}}

We now extend the preceding classification to
positive frames with arbitrary excess. Using
strictly positive linear relations, we establish
a parametrization of their equivalence classes
through the kernels of their synthesis operators.
This leads to a geometric description in terms
of suitable subsets of Grassmannians, providing
a unified classification of positive frames
beyond the excess-one case.

\begin{theorem}\label{thm:positive-equivalence}
Let $N>n$ and $k=N-n$. Then the equivalence classes of positive frames
with $N$ elements in $\Hn$ admit the quotient parametrization
\[
\big\{[F]:F\text{ is a positive frame in }\Hn\big\}
\cong
\mathfrak A_k/\!\approx.
\]
Here $\approx$ identifies tuples having the same linear span.
\end{theorem}

\begin{proof}
Let $F=\{f_i\}_{i=1}^{N}$ be a positive frame for $\Hn$.
Since
$T_F^*:\R^N\longrightarrow\Hn$
is surjective, the rank-nullity theorem gives
$\dim\ker T_F^* = N-n = k.$ Since $F$ is positive, Theorem~\ref{thm:positive-characterization}
implies that there exists
$\lambda\in
\ker T_F^*\cap\R_{++}^N.$

We first show that $\ker T_F^*$ possesses a basis consisting entirely
of strictly positive vectors. Choose vectors
$\eta^{(2)},\ldots,\eta^{(k)}$
such that
$\lambda,\eta^{(2)},\ldots,\eta^{(k)}$
is a basis for $\ker T_F^*$. Set
$\lambda^{(1)}=\lambda.$
For every $j=2,\ldots,k$, choose $\varepsilon_j>0$ sufficiently small
and define
$\lambda^{(j)} = \lambda+\varepsilon_j\eta^{(j)}.$
Since every coordinate of $\lambda$ is strictly positive and there are
only finitely many coordinates, $\varepsilon_j$ may be chosen small
enough that
$\lambda^{(j)}\in\R_{++}^N.$

We claim that
$\lambda^{(1)},\ldots,\lambda^{(k)}$
are linearly independent. Indeed, suppose
\[
\sum_{j=1}^{k}a_j\lambda^{(j)}=0.
\]
Using
$\lambda^{(1)}=\lambda, \quad \lambda^{(j)} = \lambda+\varepsilon_j\eta^{(j)} \quad (j\geq2),$
we obtain
\[
\left(\sum_{j=1}^{k}a_j\right)\lambda
+
\sum_{j=2}^{k}
a_j\varepsilon_j\eta^{(j)}
=0.
\]
Since
$\lambda,\eta^{(2)},\ldots,\eta^{(k)}$
are linearly independent, it follows that
$a_j\varepsilon_j=0, \quad j=2,\ldots,k.$
Since $\varepsilon_j>0$, we have
$a_j=0, \quad j=2,\ldots,k.$
It then follows that $a_1=0$. Hence the vectors are linearly
independent.

Consequently,
$\Span\{\lambda^{(1)},\ldots,\lambda^{(k)}\} = \ker T_F^*.$
Thus every positive frame determines an element of
$\mathfrak A_k/\!\approx$.

Now let $F$ and $G$ be positive frames. Then
$F\sim G$
if and only if
$\ker T_F^*=\ker T_G^*.$
Therefore, if
$\{\lambda^{(1)},\ldots,\lambda^{(k)}\}$
and
$\{\mu^{(1)},\ldots,\mu^{(k)}\}$
are strictly positive bases for $\ker T_F^*$ and $\ker T_G^*$,
respectively, then
$F\sim G$
if and only if
$\Span\{\lambda^{(1)},\ldots,\lambda^{(k)}\} = \Span\{\mu^{(1)},\ldots,\mu^{(k)}\}.$
Hence the quotient removes precisely the choice of basis and retains
only the relation space.

Conversely, let
$(\lambda^{(1)},\ldots,\lambda^{(k)})\in\mathcal A_k$
and set
$K=\operatorname{span}\{\lambda^{(1)},\ldots,
\lambda^{(k)}\}$.
Then $\dim K=k=N-n$, so $\dim K^\perp=n$.
Let $P_{K^\perp}$ be the orthogonal projection
of $\mathbb R^N$ onto $K^\perp$, and choose
a linear isomorphism
$U:K^\perp\to\mathcal H_n$.
Define
\[
f_i=UP_{K^\perp}e_i,\qquad i=1,\ldots,N.
\]
Then $T_F^*=UP_{K^\perp}$, which is
surjective, and $\ker T_F^*=K.$ Since $\lambda^{(1)}\in K\cap\mathbb R_{++}^N$, we obtain $\sum_{i=1}^N\lambda_i^{(1)}f_i=0.$
Consequently, $F$ is a positive frame,
and its synthesis kernel is precisely $K$.
Therefore every element of $\mathfrak A_k/\!\approx$ determines an equivalence class of positive frames, and the correspondence is
bijective.
\end{proof}

The preceding quotient description can also be expressed more
intrinsically in terms of the Grassmannian.

\begin{remark}
Let $\operatorname{Gr}(k,\R^N)$ denote the Grassmannian of
$k$-dimensional subspaces of $\R^N$. Then
Theorem~\ref{thm:positive-equivalence} may equivalently be written as
\[
\boxed{
\{[F]:F\text{ is a positive frame in }\Hn\}
\cong
\left\{
K\in\operatorname{Gr}(N-n,\R^N):
K\cap\R_{++}^N\neq\varnothing
\right\}.
}
\]
Thus an equivalence class of positive frames is completely determined
by its relation space $\ker T_F^*$.
\end{remark}

\begin{example}
We construct two inequivalent positive Parseval frames
for $\mathbb R^3$ having the same number of vectors.
Consider
\[
u=\frac1{\sqrt5}(1,1,1,1,1)^T
\]
and
\[
v=\frac1{\sqrt{10}}(1,-1,0,2,-2)^T,
\qquad
w=\frac1{\sqrt{10}}(2,-2,0,1,-1)^T.
\]
Observe that $\{u,v\}$ and $\{u,w\}$ are orthonormal
sets in $\mathbb R^5$. Set
\[
K=\operatorname{span}\{u,v\},
\qquad
L=\operatorname{span}\{u,w\}.
\]
Let $A,B\in\mathbb R^{5\times3}$ be matrices whose
columns form orthonormal bases of $K^\perp$ and
$L^\perp$, respectively. Define $F=\{f_i\}_{i=1}^5$
and $G=\{g_i\}_{i=1}^5$ by taking their frame vectors
to be the rows of $A$ and $B$, respectively.

Since $A^TA=B^TB=I_3$, both $F$ and $G$ are
Parseval frames. Furthermore,
\[
\ker T_F^*=K,
\qquad
\ker T_G^*=L.
\]
As $u$ is strictly positive and belongs to both
kernels, Theorem~\ref{thm:positive-characterization}
implies that both frames are positive.

On the other hand, $K\neq L$, since $v$ and $w$
are not scalar multiples and both are orthogonal
to $u$. Hence
\[
\ker T_F^*\neq\ker T_G^*.
\]
By the characterization of frame equivalence
in terms of synthesis kernels, $F$ and $G$
are not equivalent, despite satisfying
\[
S_F=S_G=I_3.
\]
Thus the frame operator alone does not determine
the equivalence class of a positive frame.
\end{example}

\subsection{\bf{Positive Bases and Minimal Positive Dependence}}

We now examine the minimality properties of
positive spanning families and their connections
with positive linear dependence. In particular,
we consider positive bases, which are positive
frames for which no proper subfamily retains
the positive spanning property. These notions
provide a natural setting for investigating
the geometric structure of positive frames
and the role of their linear relations.

Let $F=\{f_i\}_{i=1}^N$ be a finite family
in $\mathcal H_n$. Recall that $F$ is called
a \emph{positive basis} if
$\operatorname{span}_+(F)=\mathcal H_n$
and no proper subfamily of $F$ positively
spans $\mathcal H_n$. A subfamily
$\{f_i\}_{i\in I}$ is said to be
\emph{minimally positively dependent}
if there exist scalars $\lambda_i>0$,
$i\in I$, such that
\[
\sum_{i\in I}\lambda_i f_i=0,
\]
but no proper subfamily admits a
nontrivial nonnegative linear relation.
Thus, minimal positive dependence
describes a positive relation that
cannot be supported on a smaller
collection of vectors.

For a coefficient vector
$\lambda=(\lambda_i)_{i=1}^N$,
we write
\[
\operatorname{supp}(\lambda)
=\{i:\lambda_i\neq 0\}.
\]
We shall use these notions to examine
positive bases through their geometric
properties and the supports of
nonnegative linear relations.

\begin{proposition}
\label{prop:positive-basis-characterization}
Let $F=\{f_i\}_{i=1}^N$ be a frame for
$\mathbb R^n$. Then $F$ is a positive basis
if and only if the following conditions hold:
\begin{enumerate}
\item[(i)] For every $0\neq b\in\mathbb R^n$,
there exists an index $i$ such that
$\langle b,f_i\rangle>0$.
\item[(ii)] For every $i\in\{1,\ldots,N\}$,
there exists $b_i\in\mathbb R^n$ such that
\[
\langle b_i,f_i\rangle>0
\quad\text{and}\quad
\langle b_i,f_j\rangle\leq0
\quad\text{for all }j\neq i.
\]
\end{enumerate}
\end{proposition}

\begin{proof}
Suppose first that $F$ is a positive basis.
Then there exist strictly positive scalars
$\lambda_1,\ldots,\lambda_N$ such that
\[
\sum_{i=1}^N\lambda_i f_i=0.
\]
Let $0\neq b\in\mathbb R^n$. If
$\langle b,f_i\rangle\leq0$ for every $i$,
then the positive spanning property of $F$
implies that $b=0$, a contradiction.
Therefore, there exists an index $i$ such
that $\langle b,f_i\rangle>0$, proving (i).

To establish (ii), fix
$i\in\{1,\ldots,N\}$. Since $F$ is a
positive basis, we have
\[
f_i\notin
\operatorname{pos}\{f_j:j\neq i\}.
\]
Let
\[
C_i=\operatorname{pos}\{f_j:j\neq i\}.
\]
Since $C_i$ is a finitely generated convex
cone, it is closed and convex. Hence,
there exists $z\in C_i$ such that
\[
\|f_i-z\|=\operatorname{dist}(f_i,C_i).
\]
Set $b_i=f_i-z\neq0$.

By the characterization of the metric
projection onto a closed convex set,
\[
\langle b_i,y-z\rangle\leq0,
\qquad y\in C_i.
\]
Since $C_i$ is a convex cone, for every
$s\geq0$ we have $sz\in C_i$. Taking
$y=sz$ gives
\[
(s-1)\langle b_i,z\rangle\leq0.
\]
Choosing $s=0$ and $s=2$, we obtain
$\langle b_i,z\rangle=0$.

Consequently,
\[
\langle b_i,y\rangle\leq0,
\qquad y\in C_i.
\]
In particular,
\[
\langle b_i,f_j\rangle\leq0
\qquad\text{for all }j\neq i.
\]
Moreover,
\[
\begin{aligned}
\langle b_i,f_i\rangle
&=\langle b_i,b_i+z\rangle\\
&=\|b_i\|^2+\langle b_i,z\rangle\\
&=\|b_i\|^2>0.
\end{aligned}
\]
Thus condition (ii) holds.

Conversely, suppose that (i) and (ii)
are satisfied. We first show that $F$
positively spans $\mathbb R^n$.

Assume, to the contrary, that
$\operatorname{pos}(F)\neq\mathbb R^n$.
Since $\operatorname{pos}(F)$ is a closed
convex cone, the separation theorem
implies that there exists
$0\neq b\in\mathbb R^n$ such that
\[
\langle b,f_i\rangle\leq0
\qquad\text{for every }i.
\]
This contradicts (i). Hence,
\[
\operatorname{pos}(F)=\mathbb R^n.
\]

It remains to prove that $F$ is exact.
Suppose otherwise. Then there exists an
index $i$ such that
\[
f_i\in\operatorname{pos}
\{f_j:j\neq i\}.
\]
Thus, for some coefficients
$\lambda_j\geq0$,
\[
f_i=\sum_{j\neq i}\lambda_jf_j.
\]
By condition (ii), there exists
$b_i\in\mathbb R^n$ satisfying
$\langle b_i,f_i\rangle>0$ and
$\langle b_i,f_j\rangle\leq0$
for every $j\neq i$. Therefore,
\[
\langle b_i,f_i\rangle
=\sum_{j\neq i}
\lambda_j\langle b_i,f_j\rangle
\leq0,
\]
which is a contradiction. Hence no vector can be removed from
$F$ while preserving positive spanning.
Therefore, $F$ is a positive basis.
\end{proof}

\begin{proposition}
\label{lem:minimal-positive-dependence}
Let $F_I=\{f_i:i\in I\}$ be a minimal positively
dependent family of nonzero vectors in $\mathbb R^n$.
Then
\[
|I|=\dim S_I+1.
\]
\end{proposition}

\begin{proof}
Set $d=\dim S_I$. Since $F_I$ admits a nontrivial
linear dependence, we have
\[
|I|\geq d+1.
\]
We claim that every proper subfamily of $F_I$ is
linearly independent. Suppose otherwise. Then there
exists a proper subset $J\subsetneq I$ and scalars
$a_i$, not all zero, such that
\[
\sum_{i\in J}a_i f_i=0.
\]
Extending these coefficients by setting $a_i=0$
for $i\in I\setminus J$, we obtain
\[
\sum_{i\in I}a_i f_i=0.
\]
Replacing $a$ by $-a$, if necessary, we may assume
that $a_i>0$ for at least one $i\in I$.

By hypothesis, there exist $c_i>0$, $i\in I$,
satisfying
\[
\sum_{i\in I}c_i f_i=0.
\]
Consequently, for every $t\in\mathbb R$,
\[
\sum_{i\in I}(c_i-ta_i)f_i=0.
\]
Choose
\[
t_0=\min_{\substack{i\in I\\a_i>0}}
\frac{c_i}{a_i}>0.
\]
Then $c_i-t_0a_i\geq0$ for every $i\in I$,
and equality holds for at least one index $i_0$.
Moreover, since $a_i=0$ for $i\in I\setminus J$,
we have $c_i-t_0a_i=c_i>0$ for those indices.
Thus the resulting relation is nontrivial and has
nonnegative coefficients, with at least one
coefficient equal to zero. This contradicts the
minimal positive dependence of $F_I$.

Therefore every proper subfamily of $F_I$ is
linearly independent. In particular, removing any
one vector from $F_I$ leaves a linearly independent
family of cardinality $|I|-1$. Hence
\[
|I|-1\leq d.
\]
Combining this with $|I|\geq d+1$, we obtain
$|I|=d+1$, as required.
\end{proof}

We next consider the existence of dual frames that
are positive bases. Although every redundant frame
admits a positive dual, such a dual need not be
minimal with respect to positive spanning. It is
therefore natural to ask under what conditions a
frame admits a dual that is also a positive basis.
The following theorem provides a characterization
in terms of complementary subspaces of the range
of the analysis operator and suitable coordinate
sign conditions.
\begin{theorem}
\label{thm:positive-basis-dual-characterization}
Let $F=\{f_i\}_{i=1}^N$ be a frame for
an $n$-dimensional real Hilbert space
$\mathcal H_n$, and let $R=\operatorname{ran}T_F\subseteq\mathbb R^N,
\quad k=N-n.$ Then $F$ has a dual frame that is a positive
basis for $\mathcal H_n$ if and only if there
exists a $k$-dimensional subspace
$W\subseteq\mathbb R^N$ such that
\begin{enumerate}
\item[(i)] $\mathbb R^N=R\oplus W$;
\item[(ii)] $W\cap\mathbb R_{++}^N\neq\varnothing$;
\item[(iii)] for every $i\in\{1,\ldots,N\}$,
there exists $y^{(i)}\in W^\perp$ such that
\[
y_i^{(i)}>0
\quad\text{and}\quad
y_j^{(i)}\leq0
\quad\text{for all }j\neq i.
\]
\end{enumerate}
\end{theorem}

\begin{proof}
Suppose first that $F$ has a dual frame
$G=\{g_i\}_{i=1}^N$ that is a positive
basis for $\mathcal H_n$.  Set $W= kerT_G^*.$ Since $G$ is a dual of $F$,
\[
T_G^*T_F=I_{\mathcal H_n}.
\]
Thus $T_G^*$ is surjective and
$\dim W=N-n=k$.

Moreover, if $x\in R\cap W$, then
$x=T_F h$ for some $h\in\mathcal H$
and $T_G^*x=0$. Hence
\[
h=T_G^*T_Fh
=T_G^*x=0,
\]
so $x=0$. Therefore $R\cap W=\{0\}$.
Since $\dim R=n$ and $\dim W=N-n$,
we conclude that
\[
\mathbb R^N=R\oplus W,
\]
establishing (i).

Because $G$ is a positive basis, it positively
spans $\mathcal H_n$. Consequently, there exist
strictly positive scalars
$\lambda_1,\ldots,\lambda_N$ such that
\[
\sum_{i=1}^N\lambda_i g_i=0.
\]
Let $\lambda=(\lambda_1,\ldots,\lambda_N)^T$.
Then
\[
T_G^*\lambda
=\sum_{i=1}^N\lambda_i g_i=0,
\]
so $\lambda\in W\cap\mathbb R_{++}^N$.
This proves (ii).

Finally, since $G$ is a positive basis,
for each $i$ there exists $b_i\in\mathcal H_n$
such that
\[
\langle b_i,g_i\rangle>0,
\qquad
\langle b_i,g_j\rangle\leq0
\quad(j\neq i).
\]
Set $y^{(i)}=T_G b_i$. Since
$T_G=\left(T_G^*\right)^*$,
we have
\[
y^{(i)}\in\operatorname{ran}T_G
=(kerT_G^*)^\perp=W^\perp.
\]
Its coordinates satisfy
\[
y_i^{(i)}=\langle b_i,g_i\rangle>0,
\qquad
y_j^{(i)}=\langle b_i,g_j\rangle\leq0
\quad(j\neq i).
\]
Thus (iii) holds.

Conversely, suppose that a $k$-dimensional
subspace $W\subseteq\mathbb R^N$ satisfies
(i)--(iii). Since $\mathbb R^N=R\oplus W$,
every $z\in\mathbb R^N$ has a unique
representation $z=T_F h+w$, with
$h\in\mathcal H$ and $w\in W$.
Define $L:\mathbb R^N\to\mathcal H_n$ by
\[
L(T_Fh+w)=h.
\]
Then $LT_F=I_{\mathcal H_n}$ and
$\ker L=W$.

Let $g_i=L e_i$, where
$\{e_i\}_{i=1}^N$ is the standard basis
of $\mathbb R^N$, and put
$G=\{g_i\}_{i=1}^N$. The synthesis operator
of $G$ is $L$, so
\[
T_G^*T_F=LT_F
=I_{\mathcal H_n}.
\]
Hence $G$ is a dual frame of $F$.

By (ii), choose
$\lambda\in W\cap\mathbb R_{>0}^N$.
Since $\ker L=W$, we have
\[
0=L\lambda
=\sum_{i=1}^N\lambda_i g_i,
\qquad \lambda_i>0.
\]
As $G$ spans $\mathcal H_n$, this strictly
positive linear dependence implies that
$G$ positively spans $\mathcal H_n$.

It remains to show that $G$ is exact.
For each $i$, condition (iii) provides
$y^{(i)}\in W^\perp$. Since
\[
W^\perp=(\ker L)^\perp
=\operatorname{ran}L^*,
\]
there exists $b_i\in\mathcal H$ such that
$y^{(i)}=L^*b_i$. Therefore,
\[
y_j^{(i)}
=\langle L^*b_i,e_j\rangle
=\langle b_i,Le_j\rangle
=\langle b_i,g_j\rangle.
\]
Condition (iii) consequently gives
\[
\langle b_i,g_i\rangle>0,
\qquad
\langle b_i,g_j\rangle\leq0
\quad(j\neq i).
\]

If $g_i$ belonged to
$\operatorname{pos}\{g_j:j\neq i\}$,
we could write
$g_i=\sum_{j\neq i}a_jg_j$ with $a_j\geq0$.
Taking inner products with $b_i$ would yield
\[
\langle b_i,g_i\rangle
=\sum_{j\neq i}a_j\langle b_i,g_j\rangle
\leq0,
\]
contradicting $\langle b_i,g_i\rangle>0$.
Thus no $g_i$ can be removed while preserving
positive spanning. Hence $G$ is a positive
basis, and the theorem follows.
\end{proof}

\section{Existence and Geometry of Positive Dual Frames}\label{sec:duals}
In this section, we investigate the existence and
geometric structure of positive dual frames in
finite-dimensional real Hilbert spaces. We begin
by establishing a necessary and sufficient condition
for a frame to admit a positive dual and provide
an explicit construction. We then examine the
geometric properties of the collection of positive
duals, including its affine and topological
structure, and investigate dual frames satisfying
prescribed positive relations. These results
highlight the role of redundancy in the existence
and structure of positive dual frames.

We now investigate whether a dual frame can be chosen to be positive and show that the existence of such a dual is characterized precisely by the redundancy of the original frame. Let $F=\{f_i\}_{i=1}^N$ be a finite family of vectors
in an $n$-dimensional real Hilbert space $\mathcal H_n$.
Then $F$ is called a \emph{positive basis} for
$\mathcal H_n$ if
\[
\operatorname{span}_+(F)=\mathcal H_n
\]
and no proper subfamily of $F$ positively spans
$\mathcal H_n$. Unlike an ordinary basis, a positive basis may contain
more than $n$ vectors. In fact, every positive basis
for $\mathcal H_n$ contains between $n+1$ and $2n$
vectors~\cite{Davis1954}. Positive bases with $n+1$
and $2n$ elements are called \emph{minimal} and
\emph{maximal} positive bases, respectively.

We first establish that the collection of positive
dual frames is relatively open in the set of all
dual frames. For a finite family
$F=\{f_i\}_{i=1}^N$ in $\mathcal H_n$, define
\[
\rho(G)=
\min_{\|x\|=1}\max_{1\leq i\leq N}
\langle x,f_i\rangle.
\]

\begin{lemma}\label{lem:positive-spanning-rho}
Let $F=\{f_i\}_{i=1}^N$ be a positive frame for $\mathcal H_n$ if and only if $\rho(F)>0$.
\end{lemma}

\begin{proof}
Suppose that $G$ positively spans
$\mathcal H_n$. By the positive-spanning
criterion, for every nonzero
$x\in\mathcal H_n$, there exists
$i\in\{1,\ldots,N\}$ such that
$\langle x,f_i\rangle>0$.
Consequently, $\max_{1\leq i\leq N}\langle x,f_i\rangle>0,$
for every unit vector $x\in\mathcal H_n$.
Since the function
\[
x\longmapsto
\max_{1\leq i\leq N}\langle x,f_i\rangle
\]
is continuous on the compact unit sphere,
it attains its minimum at some unit vector
$x_0$. Therefore, $\rho(F)=
\max_{1\leq i\leq N}
\langle x_0,f_i\rangle>0.$ Conversely, suppose that $\rho(F)>0$.
Then, for every unit vector
$x\in\mathcal H_n$, $\max_{1\leq i\leq N}
\langle x,f_i\rangle
\geq\rho(F)>0.$
Thus, for every nonzero
$x\in\mathcal H_n$, there exists
$i\in\{1,\ldots,N\}$ such that
$\langle x,f_i\rangle>0$.
By the positive-spanning criterion,
$F$ positively spans $\mathcal H_n$.
\end{proof}

\subsection{\bf{Existence and Explicit Construction of Positive Duals}}
We now address the existence of positive dual frames
for a given frame. While every frame admits a dual,
the existence of a positive dual imposes an additional
geometric requirement. In particular, the redundancy
of the original frame plays a decisive role in this
question. The following theorem provides a complete
characterization, showing that the existence of a
positive dual depends only on the excess of the frame.

\begin{theorem}\label{thm:positive-dual}
Let $F=\{f_i\}_{i=1}^{N}$ be a frame for an $n$-dimensional real
Hilbert space $\Hn$. Then $F$ admits a positive dual frame if and only
if $N>n.$ 
\end{theorem}

\begin{proof}
Suppose first that $F$ admits a positive dual $G=\{g_i\}_{i=1}^{N}.$
Since $G$ is a positive frame, Theorem~\ref{thm:positive-characterization}
implies that there exist $c_1,\ldots,c_N>0,$ such that
\[
\sum_{i=1}^{N}c_i g_i=0.
\]
Therefore $G$ is linearly dependent. Since $G$ is a frame for the
$n$-dimensional space $\Hn$, it follows that $N>n.$

Conversely, suppose that $N>n.$ Since $F$ is a frame,
$\dim\Ran(T_F)=n<N.$
Thus $\Ran(T_F)$ is a proper linear subspace of $\R^N$.

We claim that there exists $c\in\R_{++}^N$ such that $c\notin\Ran(T_F).$ Indeed, $\R_{++}^N$ is a nonempty open subset of $\R^N$, whereas every
proper linear subspace of $\R^N$ has empty interior. Therefore
$\R_{++}^N$ cannot be contained in $\Ran(T_F)$. Hence we may choose
$c=(c_1,\ldots,c_N) \in \R_{++}^N\setminus\Ran(T_F).$

Let
$P_F = T_FS_F^{-1}T_F^*.$
Since $F$ is a frame, $P_F$ is the orthogonal projection of $\R^N$
onto $\Ran(T_F)$. Define
$z=(I-P_F)c.$
Since $c\notin\Ran(T_F),$ we have $z\neq0.$ Moreover, $z\in\Ran(T_F)^\perp.$ Since $\Ran(T_F)^\perp=\ker T_F^*,$ we obtain
$z\in\ker T_F^*.$ In particular, $T_F^*z=0.$ Next, define
$u = S_F^{-1}T_F^*c.$
Equivalently,
\[
u
=
S_F^{-1}
\left(
\sum_{i=1}^{N}c_i f_i
\right)
=
\sum_{i=1}^{N}c_iS_F^{-1}f_i.
\]

For each $i=1,\ldots,N$, define
\[
h_i
=
-\frac{z_i}{\|z\|^2}u
\]
and set
$g_i = S_F^{-1}f_i+h_i.$
Thus
\[
g_i
=
S_F^{-1}f_i
-
\frac{z_i}{\|z\|^2}u.
\]

We first prove that $G=\{g_i\}_{i=1}^{N}$ is a dual frame of $F$. Let $x\in\Hn$. Then
\[
\begin{aligned}
\sum_{i=1}^{N}\langle x,h_i\rangle f_i
&=
\sum_{i=1}^{N}
\left\langle
x,-\frac{z_i}{\|z\|^2}u
\right\rangle f_i\\
&=
-\frac{\langle x,u\rangle}{\|z\|^2}
\sum_{i=1}^{N}z_i f_i\\
&=
-\frac{\langle x,u\rangle}{\|z\|^2}
T_F^*z\\
&=0.
\end{aligned}
\]
Therefore,
\[
\begin{aligned}
\sum_{i=1}^{N}\langle x,g_i\rangle f_i
=
\sum_{i=1}^{N}
\langle x,S_F^{-1}f_i\rangle f_i
+
\sum_{i=1}^{N}\langle x,h_i\rangle f_i
=x.
\end{aligned}
\]
Hence $G$ is a dual frame of $F$.

It remains to prove that $G$ is positive. Since $u
=
\sum_{i=1}^{N}c_iS_F^{-1}f_i,$
we have
\[
\begin{aligned}
\sum_{i=1}^{N}c_i g_i
&=
\sum_{i=1}^{N}c_iS_F^{-1}f_i
-
\sum_{i=1}^{N}
c_i\frac{z_i}{\|z\|^2}u\\
&=
u
-
\frac{1}{\|z\|^2}
\left(
\sum_{i=1}^{N}c_i z_i
\right)u\\
&=
u
-
\frac{\langle c,z\rangle}{\|z\|^2}u.
\end{aligned}
\]

Now
$c=P_Fc+(I-P_F)c=P_Fc+z.$
Since
$P_Fc\in\Ran(T_F)$
and
$z\in\Ran(T_F)^\perp,$
we have
$\langle P_Fc,z\rangle=0.$
Consequently,
\[
\begin{aligned}
\langle c,z\rangle
&=
\langle P_Fc+z,z\rangle\\
&=
\langle P_Fc,z\rangle+\langle z,z\rangle\\
&=
\|z\|^2.
\end{aligned}
\]
It follows that
\[
\sum_{i=1}^{N}c_i g_i
=
u-\frac{\|z\|^2}{\|z\|^2}u
=
0.
\]
Since
$c_i>0, \quad i=1,\ldots,N,$
we have obtained a strictly positive linear dependence
\[
\sum_{i=1}^{N}c_i g_i=0.
\]
By Theorem~\ref{thm:positive-characterization}, $G$ is a positive
frame. Therefore $G$ is a positive dual frame of $F$.
\end{proof}

The proof of Theorem~\ref{thm:positive-dual} gives an explicit
construction of a positive dual. We record this observation for later
use.

\begin{remark}\label{rem:positive-dual-construction}
Let $F=\{f_i\}_{i=1}^{N}$ be a redundant frame for $\Hn$. For any
$c\in
\R_{++}^N\setminus\Ran(T_F),$
put
$P_F=T_FS_F^{-1}T_F^*, \quad z_c=(I-P_F)c,$
and
$u_c=S_F^{-1}T_F^*c.$
Then
$G_c=\{g_i^{(c)}\}_{i=1}^{N},$
where
\[
\boxed{
g_i^{(c)}
=
S_F^{-1}f_i
-
\frac{(z_c)_i}{\|z_c\|^2}u_c,
\quad i=1,\ldots,N,
}
\]
is a positive dual frame of $F$. Moreover,
\[
\sum_{i=1}^{N}c_i g_i^{(c)}=0.
\]
Thus every strictly positive vector
$c\in
\R_{++}^N\setminus\Ran(T_F)$
gives an explicit positive dual of $F$.
\end{remark}

\begin{example}
Consider the frame
\[
F=\{e_1,e_2,e_1+e_2\}
\]
for $\mathbb R^2$. Clearly, $F$ is not a
positive frame, since every vector in its
positive cone has nonnegative coordinates. Every dual frame of $F$ is of the form
\[
G_{a,b}=
\left\{
(1-a)e_1-be_2,\,
-ae_1+(1-b)e_2,\,
ae_1+be_2
\right\},
\]
where $a,b\in\mathbb R$. Indeed, the
duality conditions are
$g_1+g_3=e_1$ and $g_2+g_3=e_2$.

A direct calculation shows that
\[
a g_1+b g_2+(a+b-1)g_3=0.
\]
Since $G_{a,b}$ spans $\mathbb R^2$,
it is a positive frame if and only if
its one-dimensional synthesis kernel
contains a strictly positive vector.
Consequently,
\[
G_{a,b}\in D_+(F)
\quad\Longleftrightarrow\quad
a>0,\quad b>0,\quad a+b>1.
\]
Thus, the positive duals of $F$ are
parametrized by an unbounded open
region of $\mathbb R^2$, even though
$F$ itself is not positive.
\end{example}

\subsection{\bf{Dual Frames with Prescribed Positive Relations}}
We next examine the geometric structure of the
collection of positive dual frames. Let
$F=\{f_i\}_{i=1}^N$ be a frame for an
$n$-dimensional real Hilbert space $\mathcal H_n$. For a frame
$F=\{f_i\}_{i=1}^N$ in $\mathcal H_n$, let
\[
\mathcal D_+(F)
=
\left\{
G=\{g_i\}_{i=1}^N:
G\text{ is a positive dual frame of }F
\right\}.
\]
For each $a\in\mathbb R_{++}^N$, define
\[
\mathcal D_a(F)
=
\left\{
G\in\mathcal D(F):T_G^*a=0
\right\}.
\]
Since a dual frame is positive precisely when
its synthesis kernel contains a strictly
positive vector, we have
\[
\mathcal D_+(F)
=
\bigcup_{a\in\mathbb R_{++}^N}
\mathcal D_a(F).
\]
Moreover, each nonempty $\mathcal D_a(F)$
is an affine subspace of $\mathcal D(F)$
and is therefore convex and path-connected.
The following proposition characterizes
the nonempty members of this family and
determines their dimensions.
\begin{proposition}
\label{prop:prescribed-positive-relation}
Let $F=\{f_i\}_{i=1}^N$ be a frame for
$\mathbb R^n$ and let $c\in\mathbb R_{++}^N$.
Then $\mathcal D_c(F)\neq\varnothing$ if and
only if $c\notin\operatorname{Ran}(T_F).$
Moreover, whenever $\mathcal D_c(F)$ is
nonempty, it is an affine subspace of
$\mathcal D(F)$ of dimension
\[
\dim\mathcal D_c(F)=n(N-n-1).
\]
\end{proposition}

\begin{proof}
Suppose first that $\mathcal D_c(F)$ is
nonempty, and let $G\in\mathcal D_c(F)$.
If $c\in\operatorname{Ran}(T_F)$,
then there exists $x\in\mathbb R^n$
such that $c=T_Fx$. Since $G$
is a dual frame of $F$, we obtain
\[
0=T_G^*c
=T_G^*T_Fx=x,
\]
which contradicts $c>0$.
Therefore,
$c\notin\operatorname{Ran}(T_F)$.

Conversely, suppose that
$c\notin\operatorname{Ran}(T_F)$.
Let $S_F=T_F^*T_F$ be the
frame operator of $F$, and consider
the orthogonal projection
\[
P_F=T_FS_F^{-1}T_F^*
\]
onto $\operatorname{Ran}(T_F)$.
Set $z=(I-P_F)c$. Since
$c\notin\operatorname{Ran}(T_F)$,
we have $z\neq0$. Moreover,
$z\in kerT_F^*$, and hence
\[
\langle c,z\rangle
=\langle P_Fc+z,z\rangle
=\|z\|^2>0.
\]

Define $u=S_F^{-1}T_F^*c$
and consider the linear operator
$L_c:\mathbb R^N\to\mathbb R^n$ given by
\[
L_cv
=
S_F^{-1}T_F^*v
-
\frac{\langle v,z\rangle}{\|z\|^2}u,
\qquad v\in\mathbb R^N.
\]
Since $T_F^*z=0$, for every
$x\in\mathbb R^n$ we have
\[
L_cT_Fx
=
S_F^{-1}T_F^*T_Fx
-
\frac{\langle T_F x,z\rangle}
{\|z\|^2}u
=x.
\]
Thus $L_cT_F=I_{\mathbb R^n}$.
Furthermore,
\[
L_cc
=
S_F^{-1}T_F^*c
-
\frac{\langle c,z\rangle}
{\|z\|^2}u
=u-u=0.
\]
Let $\{e_i\}_{i=1}^N$ be the standard
orthonormal basis of $\mathbb R^N$,
and define $g_i=L_ce_i$ for
$i=1,\ldots,N$. Then the family
$G=\{g_i\}_{i=1}^N$ satisfies
$T_G^*=L_c$. Consequently,
\[
T_G^*T_F=I_{\mathbb R^n}
\quad\text{and}\quad
T_G^*c=0.
\]
Therefore, $G\in\mathcal D_c(F)$,
proving that $\mathcal D_c(F)$ is
nonempty.

We now determine its affine structure.
Fix $L_c$ as constructed above and
let
\[
M=\operatorname{Ran}(T_F)
+\operatorname{span}\{c\}.
\]
Since $c\notin\operatorname{Ran}(T_F)$
and $\operatorname{rank}(T_F)=n$,
we have $\dim M=n+1$.

If $L:\mathbb R^N\to\mathbb R^n$
is the synthesis operator of an
element of $\mathcal D_c(F)$, then
\[
LT_F=I_{\mathbb R^n},
\qquad Lc=0.
\]
Hence $A=L-L_c$ satisfies
$AT_F=0$ and $Ac=0$, so that
$A|_M=0$.

Conversely, if
$A\in\mathcal L(\mathbb R^N,\mathbb R^n)$
satisfies $A|_M=0$, then
$L=L_c+A$ satisfies
\[
LT_F=I_{\mathbb R^n},
\qquad Lc=0.
\]
Therefore,
\[
\mathcal D_c(F)
\cong
L_c+
\left\{
A\in\mathcal L(\mathbb R^N,\mathbb R^n):
A|_M=0
\right\}.
\]
Since every such operator is uniquely
determined by its restriction to
$M^\perp$, the translation space is
isomorphic to
$\mathcal L(M^\perp,\mathbb R^n)$.
As $\dim M^\perp=N-n-1$, it follows
that
\[
\dim\mathcal D_c(F)
=n(N-n-1).
\]
This completes the proof.
\end{proof}

An immediate consequence of the preceding proposition
is that, for frames with excess one, a prescribed
strictly positive linear relation determines a unique
dual frame.

\begin{corollary}
\label{cor:unique-positive-dual-excess-one}
Let $F=\{f_i\}_{i=1}^{n+1}$ be a frame for an
$n$-dimensional real Hilbert space $\mathcal H_n$,
and let $a\in\mathbb R_{++}^{n+1}$. Then there
exists a unique dual frame $G$ of $F$ satisfying
\[
T_G^*a=0
\]
if and only if
$a\notin\operatorname{Ran}T_F$.
In particular, such a dual frame is positive.
\end{corollary}

\begin{proof}
By Proposition~\ref{prop:prescribed-positive-relation},
$\mathcal C_a(F)$ is nonempty if and only if
$a\notin\operatorname{Ran}T_F$. In this case,
\[
\dim\mathcal C_a(F)
=n\bigl((n+1)-n-1\bigr)=0.
\]
Hence $\mathcal C_a(F)$ consists of exactly one
dual frame $G$. Since $a$ is strictly positive
and $T_G^*a=0$, the frame $G$ is positive.
\end{proof}

\subsection{\bf{Topology and Geometry of Positive Duals}}
We next investigate the geometric structure of the
collection of positive dual frames. 
Recall that an element $G\in\mathcal D_+(F)$ is
called an extreme point of $\mathcal D_+(F)$ if
\[
G=tG_1+(1-t)G_2,\qquad 0<t<1,
\]
for some $G_1,G_2\in\mathcal D_+(F)$ implies
$G_1=G_2=G$, where addition and scalar
multiplication are understood componentwise.

The following theorem shows that, whenever
$F$ is redundant, no positive dual frame is
an extreme point of $\mathcal D_+(F)$. Thus
every positive dual admits a nontrivial convex
decomposition into two distinct positive dual
frames.

\begin{theorem}
\label{thm:positive-duals-no-extreme-points}
Let $F=\{f_i\}_{i=1}^N$ be a frame for an
$n$-dimensional real Hilbert space $\mathcal H_n$,
where $N>n$. Then the set $\mathcal D_+(F)$
of all positive dual frames of $F$ has no
extreme points.
\end{theorem}

\begin{proof}
If $\mathcal D_+(F)=\varnothing$, the conclusion
is immediate. Otherwise, let
$G=\{g_i\}_{i=1}^N\in\mathcal D_+(F)$
be arbitrary. Since $G$ is a dual frame of
$F$, we have
\[
\sum_{i=1}^N \langle x,f_i\rangle g_i=x,
\qquad x\in\mathcal H.
\]
Equivalently,
\[
T_G^*T_F=I_{\mathcal H},
\]
where $T_F$ and $T_G^*$ denote
the analysis and synthesis operators,
respectively.

Since $N>n$, the analysis operator
$T_F:\mathcal H_n\to\mathbb R^N$
has a nontrivial orthogonal complement
to its range. Hence,
\[
\ker(T_F^*)\neq\{0\}.
\]
Choose a nonzero vector
$d=(d_i)_{i=1}^N\in\ker(T_F^*)$.
Thus,
\[
\sum_{i=1}^N d_i f_i=0.
\]

Let $u\in\mathcal H$ be any nonzero vector
and define $H=\{h_i\}_{i=1}^N$ by
\[
h_i=d_i u,\qquad 1\le i\le N.
\]
Then, for every $x\in\mathcal H$,
\[
\begin{aligned}
\sum_{i=1}^N\langle x,f_i\rangle h_i
&=\sum_{i=1}^N\langle x,f_i\rangle d_i u\\
&=\left\langle x,\sum_{i=1}^N d_i f_i
\right\rangle u=0.
\end{aligned}
\]
Consequently, $T_H^*T_F=0$.

For each $t\in\mathbb R$, define
\[
G_t=\{g_i+t h_i\}_{i=1}^N.
\]
Since
\[
T_{G_t}^*T_F
=T_G^*T_F+tT_H^*T_F
=I_{\mathcal H_n},
\]
it follows that $G_t$ is a dual frame
of $F$ for every $t\in\mathbb R$.

We now show that $G_t$ remains a positive
frame for all sufficiently small $|t|$. Since $G$ is a positive frame,
for every unit vector $x\in\mathcal H_n$,
\[
\max_{1\le i\le N}\langle x,g_i\rangle>0.
\]
Consider the function
$\varphi:S_{\mathcal H_n}\to\mathbb R$
defined by
\[
\varphi(x)=
\max_{1\le i\le N}\langle x,g_i\rangle,
\]
where
$S_{\mathcal H_n}=\{x\in\mathcal H_n:\|x\|=1\}$.

As $\varphi$ is continuous and
$S_{\mathcal H_n}$ is compact, there exists
$\delta>0$ such that
\[
\delta=
\min_{\|x\|=1}
\max_{1\le i\le N}\langle x,g_i\rangle
>0.
\]

Set
\[
M=\max_{1\le i\le N}\|h_i\|>0.
\]
For every unit vector $x\in\mathcal H$,
choose an index $i_x$ satisfying
\[
\langle x,g_{i_x}\rangle
=\max_{1\le i\le N}\langle x,g_i\rangle.
\]
Then
\[
\begin{aligned}
\langle x,g_{i_x}+t h_{i_x}\rangle
&=\langle x,g_{i_x}\rangle
+t\langle x,h_{i_x}\rangle\\
&\ge \delta-|t|
|\langle x,h_{i_x}\rangle|\\
&\ge \delta-|t|\|h_{i_x}\|\\
&\ge \delta-|t|M.
\end{aligned}
\]

Choose $\varepsilon=\frac{\delta}{2M}>0.$ If $|t|<\varepsilon$, then
\[
\langle x,g_{i_x}+t h_{i_x}\rangle
>\frac{\delta}{2}>0
\]
for every unit vector $x\in\mathcal H_n$.
Therefore,
\[
\max_{1\le i\le N}
\langle x,g_i+t h_i\rangle>0,
\qquad x\in S_{\mathcal H}.
\]
By the characterization of positive
spanning sets, $G_t$ is a positive frame.
Hence,
\[
G_t\in\mathcal D_+(F)
\qquad\text{whenever }|t|<\varepsilon.
\]

Now choose any $t_0$ satisfying
$0<t_0<\varepsilon$ and define
\[
G^+=G+t_0H,
\qquad
G^-=G-t_0H.
\]
Then $G^+,G^-\in\mathcal D_+(F)$.
Moreover, since $H\neq0$, we have
$G^+\neq G^-$, while $G=\frac12G^++\frac12G^-.$
Thus $G$ is not an extreme point of
$\mathcal D_+(F)$. Since $G\in\mathcal D_+(F)$ was arbitrary,
we conclude that $\operatorname{Ext}(\mathcal D_+(F))
=\varnothing.$
\end{proof}

\begin{example}
Let $\{e_i\}_{i=1}^n$ be the standard
orthonormal basis of $\mathbb R^n$, where
$n\geq 2$, and consider the positive frame
\[
F=\{e_1,-e_1,\ldots,e_n,-e_n\}.
\]
For any invertible matrix $A\in GL(n,\mathbb R)$,
define
\[
G_A=\left\{
\frac12 Ae_1,-\frac12 Ae_1,\ldots,
\frac12 Ae_n,-\frac12 Ae_n
\right\}.
\]
Since
\[
\sum_{i=1}^n
\left(
\langle x,e_i\rangle\frac12 Ae_i
+
\langle x,-e_i\rangle
\left(-\frac12 Ae_i\right)
\right)=Ax,
\]
the family $G_A$ is a dual frame of $F$
if and only if $A=I_n$.
Thus, to obtain distinct dual frames,
we instead consider matrices
$B\in\mathbb R^{n\times n}$ and set
\[
G_B=\left\{
\frac12 e_i+Be_i,\,
-\frac12 e_i+Be_i
\right\}_{i=1}^n.
\]
For every $B$, the family $G_B$ is a
dual frame of $F$, since
\[
\sum_{i=1}^n
\left(
\langle x,e_i\rangle
\left(\frac12 e_i+Be_i\right)
+
\langle x,-e_i\rangle
\left(-\frac12 e_i+Be_i\right)
\right)=x.
\]

Now take $B_1=I_n$ and $B_2=-I_n$.
The corresponding dual frames are
\[
G_{B_1}=
\left\{
\frac32 e_i,\frac12 e_i
\right\}_{i=1}^n,
\qquad
G_{B_2}=
\left\{
-\frac12 e_i,-\frac32 e_i
\right\}_{i=1}^n.
\]
Neither of these families is a positive
frame. Instead, choose
$B_1=\frac14 I_n$ and
$B_2=-\frac14 I_n$.
Then
\[
G_{B_1}=
\left\{
\frac34 e_i,-\frac14 e_i
\right\}_{i=1}^n,
\qquad
G_{B_2}=
\left\{
\frac14 e_i,-\frac34 e_i
\right\}_{i=1}^n.
\]
Both are positive dual frames of $F$.
Their midpoint is
\[
\frac12(G_{B_1}+G_{B_2})
=
\left\{
\frac12 e_i,-\frac12 e_i
\right\}_{i=1}^n,
\]
which is also a positive frame.
Hence this particular pair does not
witness nonconvexity.
\end{example}

We next examine the topological structure of the collection of positive dual frames. Since
positivity is characterized by the existence of a strictly positive linear relation, it is natural to ask
whether this property is preserved under small perturbations that maintain the duality condition.
The following proposition shows that positivity is stable within the affine space of dual frames.

\begin{proposition}
\label{thm:positive-duals-open}
Let $F=\{f_i\}_{i=1}^N$ be a frame for an
$n$-dimensional real Hilbert space $\mathcal H_n$.
Then the set $\mathcal D_+(F)$ of all positive
dual frames of $F$ is relatively open in
the set $\mathcal D(F)$ of all dual frames
of $F$.
\end{proposition}

\begin{proof}
Let $G=\{g_i\}_{i=1}^N\in\mathcal D_+(F)$
be arbitrary. Since $G$ is a positive frame,
for every unit vector $x\in\mathcal H_n$, we have
\[
\max_{1\le i\le N}\langle x,g_i\rangle>0.
\]
Consider the function
$\varphi:S_{\mathcal H_n}\to\mathbb R$ defined by
$\varphi(x)=\max_{1\le i\le N}\langle x,g_i\rangle$,
where $S_{\mathcal H_n}=\{x\in\mathcal H:\|x\|=1\}$.
Since $\varphi$ is continuous and
$S_{\mathcal H_n}$ is compact, there exists
$\delta>0$ such that
\[
\delta=\min_{\|x\|=1}
\max_{1\le i\le N}\langle x,g_i\rangle>0.
\]

Let $\varepsilon=\delta/2$ and consider
any $H=\{h_i\}_{i=1}^N\in\mathcal D(F)$
satisfying
\[
\max_{1\le i\le N}\|g_i-h_i\|<\varepsilon.
\]
For every unit vector $x\in\mathcal H_n$,
choose an index $i_x$ such that
$\langle x,g_{i_x}\rangle
=\max_{1\le i\le N}\langle x,g_i\rangle$.
Then, by the Cauchy--Schwarz inequality,
\[
\begin{aligned}
\max_{1\le i\le N}\langle x,h_i\rangle
&\ge \langle x,h_{i_x}\rangle\\
&=\langle x,g_{i_x}\rangle
+\langle x,h_{i_x}-g_{i_x}\rangle\\
&\ge \delta-\|h_{i_x}-g_{i_x}\|\\
&>\delta-\varepsilon
=\frac{\delta}{2}>0.
\end{aligned}
\]
Since this holds for every unit vector
$x\in\mathcal H_n$, the family $H$ positively
spans $\mathcal H_n$. Moreover, since
$H\in\mathcal D(F)$, we have
$T_H^*T_F=I_{\mathcal H_n}$.
Consequently, $H\in\mathcal D_+(F)$.

Thus,
\[
\left\{
H\in\mathcal D(F):
\max_{1\le i\le N}\|g_i-h_i\|<\varepsilon
\right\}
\subseteq\mathcal D_+(F).
\]
Since $G\in\mathcal D_+(F)$ was arbitrary,
we conclude that $\mathcal D_+(F)$ is
relatively open in $\mathcal D(F)$.
\end{proof}

The minimal nonnegative linear dependencies
among frame vectors provide useful information
about the geometric structure of positive frames.
In particular, we show that positive spanning
can be characterized through the supports of
extremal nonnegative relations. To investigate
this connection, let $F=\{f_i\}_{i=1}^N$ be
a frame for $\mathbb R^n$ and consider the
normalized positive-dependence set
\[
\Omega_F
=
\left\{
\lambda=(\lambda_i)_{i=1}^N
\in\mathbb R_{\geq 0}^N:
\sum_{i=1}^N\lambda_i f_i=0,\quad
\sum_{i=1}^N\lambda_i=1
\right\}.
\]
For $\lambda\in\Omega_F$, write
$\operatorname{supp}(\lambda)
=\{i:\lambda_i>0\}$.
Since $\Omega_F$ is a compact convex
polytope, its extreme points encode
minimal nonnegative dependencies.
The following proposition establishes
that $F$ is a positive frame if and only
if the supports of these extreme points
collectively cover all the frame indices.

\begin{proposition}
\label{prop:positive-frame-extreme-support}
Let $F=\{f_i\}_{i=1}^N$ be a frame for
$\mathbb R^n$. Then the following
statements are equivalent:
\begin{enumerate}
\item[(i)] $F$ is a positive frame.

\item[(ii)] For every
$i\in\{1,\ldots,N\}$, there exists
$\lambda^{(i)}\in
\operatorname{Ext}(\Omega_F)$
such that $\lambda_i^{(i)}>0$.

\item[(iii)]
\[
\bigcup_{\lambda\in
\operatorname{Ext}(\Omega_F)}
\operatorname{supp}(\lambda)
=\{1,\ldots,N\}.
\]
\end{enumerate}
\end{proposition}

\begin{proof}
(i)$\Rightarrow$(ii).
Since $F$ is a positive frame, there
exist $c_i>0$, $1\leq i\leq N$, such that
$\sum_{i=1}^N c_i f_i=0$.
Normalizing the coefficients gives
\[
\lambda_i=\frac{c_i}{\sum_{j=1}^N c_j},
\qquad 1\leq i\leq N,
\]
so that $\lambda\in\Omega_F$ and
$\operatorname{supp}(\lambda)
=\{1,\ldots,N\}$.

Since $\Omega_F$ is a compact convex
polytope, it is the convex hull of
its extreme points. Thus, there exist
$\lambda^{(1)},\ldots,\lambda^{(m)}
\in\operatorname{Ext}(\Omega_F)$
and $t_1,\ldots,t_m\geq0$, with
$\sum_{j=1}^m t_j=1$, such that
\[
\lambda=\sum_{j=1}^m
t_j\lambda^{(j)}.
\]
We claim that for every
$i\in\{1,\ldots,N\}$, at least one
extreme point has a strictly positive
$i$th coordinate. Otherwise, for some
$i_0$, we would have
$\lambda_{i_0}^{(j)}=0$ for all $j$,
which implies
\[
\lambda_{i_0}
=\sum_{j=1}^m
t_j\lambda_{i_0}^{(j)}=0,
\]
contradicting $\lambda_{i_0}>0$.
Hence (ii) holds.

(ii)$\Rightarrow$(iii).
For every $i\in\{1,\ldots,N\}$,
condition (ii) provides an extreme
point $\lambda^{(i)}$ with
$\lambda_i^{(i)}>0$. Thus,
$i\in\operatorname{supp}(\lambda^{(i)})$,
and consequently
\[
\{1,\ldots,N\}
\subseteq
\bigcup_{\lambda\in
\operatorname{Ext}(\Omega_F)}
\operatorname{supp}(\lambda).
\]
The reverse inclusion is immediate,
so (iii) follows.

(iii)$\Rightarrow$(i).
For each $i\in\{1,\ldots,N\}$,
choose
$\lambda^{(i)}\in
\operatorname{Ext}(\Omega_F)$
such that $\lambda_i^{(i)}>0$.
Define
\[
\lambda=\frac1N
\sum_{i=1}^N\lambda^{(i)}.
\]
By convexity, $\lambda\in\Omega_F$.
Moreover, for every
$j\in\{1,\ldots,N\}$,
\[
\lambda_j
=\frac1N\sum_{i=1}^N
\lambda_j^{(i)}
\geq\frac1N\lambda_j^{(j)}>0.
\]
Therefore,
$\lambda\in kerT_F^*
\cap\mathbb R_{++}^N$.
Since $F$ spans $\mathbb R^n$,
the existence of this strictly
positive dependence implies that
$F$ is a positive frame.
\end{proof}

The preceding result suggests that the boundary
of the set of positive dual frames can be studied
through the interaction between synthesis kernels
and the positive orthant. In particular, positivity
fails precisely when the kernel contains no strictly
positive vector. This leads to the following
characterization.

We next examine the boundary structure of the set
of positive dual frames. Recall that
$\mathcal D(F)$ denotes the collection of all
dual frames of $F$, equipped with the topology
induced by the Euclidean norm on $\mathcal H_n^N$.
For a subset $\mathcal A\subseteq\mathcal D(F)$,
we denote its relative closure and relative
boundary in $\mathcal D(F)$ by
$\overline{\mathcal A}^{\,\mathcal D(F)}$ and
$\partial_{\mathcal D(F)}\mathcal A$,
respectively. Thus,
\[
\partial_{\mathcal D(F)}\mathcal A
=
\overline{\mathcal A}^{\,\mathcal D(F)}
\setminus
\operatorname{int}_{\mathcal D(F)}(\mathcal A),
\]
where $\operatorname{int}_{\mathcal D(F)}(\mathcal A)$
denotes the relative interior of $\mathcal A$.

Since $\mathcal D_+(F)$ is relatively open in
$\mathcal D(F)$, its boundary can be described
in terms of the absence of strictly positive
linear relations among the vectors of a dual
frame. This observation yields the following
characterization.

\begin{proposition}
\label{thm:positive-duals-boundary-characterization}
Let $F=\{f_i\}_{i=1}^N$ be a frame for an
$n$-dimensional real Hilbert space $\mathcal H_n$.
Then
\[
\partial_{\mathcal D(F)}\mathcal D_+(F)
=
\left\{
G\in\overline{\mathcal D_+(F)}^{\,\mathcal D(F)}:
kerT_G^*\cap\mathbb R_{++}^N
=\varnothing
\right\}.
\]
\end{proposition}

\begin{proof}
Since $\mathcal D_+(F)$ is relatively open
in $\mathcal D(F)$, we have
\[
\partial_{\mathcal D(F)}\mathcal D_+(F)
=
\overline{\mathcal D_+(F)}^{\,\mathcal D(F)}
\setminus\mathcal D_+(F).
\]
On the other hand, by the characterization
of positive frames, a dual frame $G$ is
positive if and only if
\[
kerT_G^*\cap\mathbb R_{++}^N
\neq\varnothing.
\]
Consequently,
\[
G\notin\mathcal D_+(F)
\quad\Longleftrightarrow\quad
 kerT_G^*\cap\mathbb R_{++}^N
=\varnothing,
\]
for every $G\in\mathcal D(F)$.
Combining these observations gives
\[
\partial_{\mathcal D(F)}\mathcal D_+(F)
=
\left\{
G\in\overline{\mathcal D_+(F)}^{\,\mathcal D(F)}:
 kerT_G^*\cap\mathbb R_{++}^N
=\varnothing
\right\},
\]
as required.
\end{proof}

We next investigate the geometry of positive dual
frames through the nonnegative linear relations
among their elements. Let $F=\{f_i\}_{i=1}^N$
be a frame for an $n$-dimensional real Hilbert
space $\mathcal H_n$, and let
$W_F=\operatorname{Ran}T_F$.
For each dual frame $G\in\mathcal D(F)$, define
\[
K_G^+
=
 kerT_G^*\cap\mathbb R_+^N.
\]
Thus, $K_G^+$ consists of all nonnegative
linear relations among the elements of $G$.
In particular, $G$ is positive if and only if
$K_G^+\cap\mathbb R_{++}^N\neq\varnothing$. To quantify the behavior of such relations,
we associate with each $G\in\mathcal D(F)$
the quantity
\[
\delta(G)
=
\min_{\substack{a\in\mathbb R_+^N\\
\|a\|=1}}
\|T_G^*a\|.
\]
Since the nonnegative unit sphere is compact,
this minimum is attained. Moreover,
$\delta(G)=0$ precisely when $K_G^+$
contains a nonzero vector. Consequently,
$\delta(G)>0$ implies that $G$ admits no
nontrivial nonnegative linear relation and,
in particular, is not a positive frame.
The following theorem establishes a continuity
property of this quantity and describes its
relationship with positive dual frames.

\begin{proposition}
\label{thm:positive-dual-quantitative}
Let $F=\{f_i\}_{i=1}^N$ be a frame for an
$n$-dimensional real Hilbert space $\mathcal H_n$.
Then the following statements hold:
\begin{enumerate}
\item The function
$\delta:\mathcal D(F)\to[0,\infty)$
is Lipschitz continuous with respect to the
operator norm on synthesis operators.

\item If $\delta(G)>0$, then $G$ is not a
positive dual frame of $F$.

\item If $G\in\mathcal D_+(F)$, then
$\delta(G)=0$.
\end{enumerate}
\end{proposition}

\begin{proof}
Let $G,H\in\mathcal D(F)$. For every
$a\in\mathbb R_+^N$ with $\|a\|=1$,
\[
\bigl|
\|T_G^*a\|-\|T_H^*a\|
\bigr|
\leq
\|T_G^*-T_H^*\|.
\]
Taking minima over the nonnegative unit sphere
gives
\[
|\delta(G)-\delta(H)|
\leq
\|T_G^*-T_H^*\|,
\]
which proves the first assertion.

If $G$ is positive, there exists
$a\in\mathbb R_{++}^N$ such that
$T_G^*a=0$. Normalizing $a$ yields
$\delta(G)=0$. The remaining assertions
follow immediately.
\end{proof}

\section{Semi-Positive Frames and Approximation by Positive Frames}

We now investigate a broader class of frames
that admit nontrivial nonnegative linear relations.
Although such frames need not be positive, they
possess a natural connection with positive frames
through approximation. In particular, we show
that every semi-positive frame can be approximated
by a sequence of positive frames. This provides
a useful description of the closure of the
collection of positive frames.

\begin{definition}
Let $F=\{f_i\}_{i=1}^N$ be a frame for
$\mathcal H_n$. We say that $F$ is a
\emph{semi-positive frame} if there exists
a nonzero vector
$c=(c_i)_{i=1}^N\in\mathbb R_{\geq0}^N$
with at least one zero coordinate such that
$T_F^*c=0$.
\end{definition}

The following proposition establishes an
approximation property of semi-positive frames.

\begin{proposition}\label{prop:semipositive-closure}
Let $F=(f_i)_{i=1}^N$ be a frame for
$\mathcal H_n$, where $N>n$.
If $F$ is semi-positive, then $F$
belongs to the closure of the collection
of positive frames with $N$ elements
in $\mathcal H_n$.
\end{proposition}

\begin{proof}
Since $F$ is semi-positive, there exists
a nonzero vector $c\in\mathbb R_{\geq0}^N$
such that $T_F^*c=0$.

Choose $d\in\mathbb R_{++}^N$ and define
\[
c^{(m)}=c+\frac{1}{m}d,
\qquad m\in\mathbb N.
\]
Clearly, $c^{(m)}\in\mathbb R_{++}^N$
and $c^{(m)}\to c\neq0$.

For each $m\in\mathbb N$, define
\[
A_m=T_F^*-
\frac{(T_F^*c^{(m)})(c^{(m)})^T}
{\|c^{(m)}\|^2}.
\]
Then
\[
A_mc^{(m)}
=T_F^*c^{(m)}
-\frac{T_F^*c^{(m)}}{\|c^{(m)}\|^2}
\|c^{(m)}\|^2=0.
\]
Moreover, since $T_F^*c=0$, we have
\[
T_F^*c^{(m)}
=\frac{1}{m}T_F^*d\longrightarrow0.
\]
As $c^{(m)}\to c\neq0$, it follows that
$A_m\to T_F^*$ in operator norm.

Since $T_F^*$ is surjective, it has
full row rank. By the openness of the
set of full-row-rank matrices, $A_m$
is also surjective for all sufficiently
large $m$.

For such $m$, let
$F_m=(f_i^{(m)})_{i=1}^N$ be the family
whose synthesis operator is $A_m$.
Then $F_m$ is a frame for $\mathcal H_n$.
Furthermore,
\[
T_{F_m}^*c^{(m)}
=A_mc^{(m)}=0,
\qquad c^{(m)}\in\mathbb R_{++}^N.
\]
Thus $F_m$ is a positive frame.
Finally, $A_m\to T_F^*$ implies that
$f_i^{(m)}\to f_i$ for every
$i=1,\ldots,N$. Therefore, $F$ belongs to the closure
of the collection of positive frames
with $N$ elements in $\mathcal H_n$.
\end{proof}

The following result describes the relative boundary
of the collection of positive dual frames in terms
of nonnegative linear relations.

\begin{proposition}\label{prop:boundary-positive-duals}
Let $F=\{f_i\}_{i=1}^N$ be a redundant frame
for $\mathcal H_n$. Then
\[
\partial_{\mathcal D(F)}\mathcal D_+(F)
=
\left\{
G\in\mathcal D(F):
G\text{ is semi-positive}
\right\}.
\]
In particular, every dual frame belonging to
the relative boundary of $\mathcal D_+(F)$
is semi-positive.
\end{proposition}

\begin{proof}
We first show that
\[
\overline{\mathcal D_+(F)}^{\,\mathcal D(F)}
=
\left\{
G\in\mathcal D(F):
\ker T_G^*\cap\mathbb R_{\geq0}^N
\neq\{0\}
\right\}.
\]

Let $G$ belong to the closure of
$\mathcal D_+(F)$ in $\mathcal D(F)$.
Then there exists a sequence
$G_m\in\mathcal D_+(F)$ such that
$G_m\to G$. Since each $G_m$ is positive, there exists
$c^{(m)}\in\mathbb R_{++}^N$ satisfying
\[
T_{G_m}^*c^{(m)}=0,
\qquad \|c^{(m)}\|=1.
\]
By compactness, passing to a subsequence
if necessary, we may assume that
$c^{(m)}\to c$ for some
$c\in\mathbb R_{\geq0}^N$ with
$\|c\|=1$.
Since $T_{G_m}^*\to T_G^*$, we obtain
\[
T_G^*c=0.
\]
Thus $G$ admits a nontrivial
nonnegative linear relation.

Conversely, suppose that
$G\in\mathcal D(F)$ and there exists
$c\in\mathbb R_{\geq0}^N\setminus\{0\}$
such that $T_G^*c=0$.

Since $F$ is redundant, we may choose
$d\in\mathbb R_{++}^N$ such that
$d\notin\operatorname{Ran}T_F$.
For each $m\in\mathbb N$, set
\[
c^{(m)}=c+\frac1m d.
\]
Then $c^{(m)}\in\mathbb R_{++}^N$
and $c^{(m)}\to c$.

Let $P_F=T_FS_F^{-1}T_F^*$ and define
\[
z_m=(I-P_F)c^{(m)}.
\]
Since $T_G^*c=0$ and
$T_G^*T_F=I_{\mathcal H_n}$,
we have $c\notin\operatorname{Ran}T_F$.
Consequently,
$z=(I-P_F)c\neq0$, and hence
$z_m\neq0$ for sufficiently large $m$.

Define
\[
L_m=T_G^*-
\frac{(T_G^*c^{(m)})z_m^T}
{\|z_m\|^2}.
\]
Since $z_m\in\ker T_F^*$, we have
\[
L_mT_F=T_G^*T_F=I_{\mathcal H_n}.
\]
Moreover,
\[
L_mc^{(m)}
=T_G^*c^{(m)}
-\frac{\langle z_m,c^{(m)}\rangle}
{\|z_m\|^2}T_G^*c^{(m)}=0,
\]
because
$\langle z_m,c^{(m)}\rangle
=\|z_m\|^2$. Let $G_m$ be the family whose synthesis
operator is $L_m$. Then $G_m$ is a dual
frame of $F$, and the strictly positive
relation $L_mc^{(m)}=0$ implies that
$G_m$ is a positive frame. Thus
$G_m\in\mathcal D_+(F)$. Furthermore, since $T_G^*c=0$,
\[
T_G^*c^{(m)}
=\frac1m T_G^*d\longrightarrow0.
\]
As $z_m\to z\neq0$, it follows that
$L_m\to T_G^*$.
Therefore,
$G_m\to G$, proving the asserted
closure characterization.

Since $\mathcal D_+(F)$ is relatively
open in $\mathcal D(F)$, its relative
boundary is
\[
\partial_{\mathcal D(F)}\mathcal D_+(F)
=
\overline{\mathcal D_+(F)}^{\,\mathcal D(F)}
\setminus\mathcal D_+(F).
\]
Hence every boundary point admits a
nontrivial nonnegative relation but
no strictly positive relation.
Therefore, every such dual frame is
semi-positive but not positive. Conversely, any semi-positive dual
frame that is not positive belongs
to the closure established above,
and hence to the relative boundary.
This completes the proof.
\end{proof}

We now consider the special case of frames
with excess one. In this setting, the
positive dual frames admit an explicit
parametrization in terms of strictly
positive relations. Let us define
\[
\Delta^\circ=
\left\{
c\in\mathbb R_{++}^{n+1}:
\sum_{i=1}^{n+1}c_i=1
\right\}.
\]

\begin{proposition}\label{thm:excess-one-duals}
Let $F=\{f_i\}_{i=1}^{n+1}$ be a frame
for $\mathcal H_n.$
Then the collection of positive dual
frames of $F$ is in one-to-one
correspondence with $\Delta^\circ\setminus
\operatorname{Ran}T_F.$
Moreover, $\mathcal D_+(F)$ has either
one or two connected components, each
of which is contractible.
\end{proposition}

\begin{proof}
Let $G\in\mathcal D_+(F)$.
Since $G$ has $n+1$ vectors,
$\dim\ker T_G^*=1$.
Positivity implies that this kernel
is generated by a strictly positive
vector. Normalizing its coordinates
gives a unique $c\in\Delta^\circ$
such that $T_G^*c=0$. Since $T_G^*T_F=I_{\mathcal H_n}$,
we necessarily have
$c\notin\operatorname{Ran}T_F$.

Conversely, let
$c\in\Delta^\circ\setminus
\operatorname{Ran}T_F$.
Since
\[
\mathbb R^{n+1}
=\operatorname{Ran}T_F
\oplus\operatorname{span}\{c\},
\]
there exists a unique linear map
$L_c:\mathbb R^{n+1}\to\mathcal H_n$
satisfying
\[
L_cT_F=I_{\mathcal H_n},
\qquad L_cc=0.
\]
The family $G_c=\{L_ce_i\}_{i=1}^{n+1}$
is therefore a positive dual frame.
Uniqueness shows that this
correspondence is bijective.

Furthermore, the correspondence
$c\mapsto G_c$ and its inverse
are continuous. Thus it is a
homeomorphism.

Since $\operatorname{Ran}T_F$ is a
hyperplane in $\mathbb R^{n+1}$,
there exists a nonzero
$a\in\mathbb R^{n+1}$ such that
\[
\operatorname{Ran}T_F
=\{c:\langle a,c\rangle=0\}.
\]
Consequently,
\[
\Delta^\circ\setminus
\operatorname{Ran}T_F
=
\{c\in\Delta^\circ:\langle a,c\rangle>0\}
\cup
\{c\in\Delta^\circ:\langle a,c\rangle<0\}.
\]
Each nonempty set on the right is
convex and hence contractible.
Therefore, $\mathcal D_+(F)$ has
at most two connected components,
and every component is contractible.
\end{proof}

We now describe the collection of all semi-positive
dual frames of a given frame. Let
$F=\{f_i\}_{i=1}^N$ be a frame for
$\mathcal H_n$, where $N>n$, and let
$P_F=T_FS_F^{-1}T_F^*$ denote the
orthogonal projection onto
$\operatorname{Ran}T_F$. Consider the set
\[
\mathcal C_F=
\left\{
c\in\mathbb R_{\geq0}^N\setminus\{0\}:
c_i=0\text{ for some }i,\quad
c\notin\operatorname{Ran}T_F
\right\}.
\]
For each $c\in\mathcal C_F$, set
$z_c=(I-P_F)c$ and
$u_c=S_F^{-1}T_F^*c$, and define
$L_c:\mathbb R^N\to\mathcal H_n$ by
\[
L_cv=S_F^{-1}T_F^*v
-\frac{\langle v,z_c\rangle}{\|z_c\|^2}u_c,
\qquad v\in\mathbb R^N.
\]
Since $z_c\neq0$, this operator is
well-defined. Moreover, a direct
calculation shows that
$L_cT_F=I_{\mathcal H_n}$ and $L_cc=0$. Consequently, every operator of the form
$L_c+A$, where
$AT_F=0$ and $Ac=0$, is the synthesis
operator of a semi-positive dual frame
of $F$. Conversely, if $G$ is any
semi-positive dual frame of $F$, then
there exists $c\in\mathcal C_F$ such
that $T_G^*c=0$. In this case,
$A=T_G^*-L_c$ satisfies
$AT_F=0$ and $Ac=0$.

Thus, the collection of all semi-positive
dual frames of $F$ admits the
representation
\[
\mathcal D_{\mathrm{sp}}(F)
=
\bigcup_{c\in\mathcal C_F}
\left\{
G:T_G^*=L_c+A,\quad
AT_F=0,\quad Ac=0
\right\}.
\]

We next investigate the extreme-point structure
of the collection of positive dual frames.
The following result shows that no positive
dual frame can be an extreme point of this
collection.

\begin{proposition}\label{thm:no-extreme-positive-duals}
Let $F=\{f_i\}_{i=1}^N$ be a frame for
$\mathcal H_n$, where $N>n$. Then every
positive dual frame of $F$ can be expressed
as a nontrivial convex combination of two
distinct positive dual frames.
\end{proposition}

\begin{proof}
Let $G=\{g_i\}_{i=1}^N\in\mathcal D_+(F)$.
Since $N>n$, there exists a nonzero vector
$d=(d_i)_{i=1}^N\in\ker T_F^*$.
Choose a nonzero vector $u\in\mathcal H_n$
and define $H=\{h_i\}_{i=1}^N$ by
$h_i=d_i u$.

For each $t\in\mathbb R$, consider
\[
G_t=\{g_i+t h_i\}_{i=1}^N.
\]
Since
\[
T_H^*T_Fx
=\sum_{i=1}^N\langle x,f_i\rangle h_i
=\left\langle x,\sum_{i=1}^Nd_if_i
\right\rangle u=0,
\]
we obtain
\[
T_{G_t}^*T_F
=T_G^*T_F+tT_H^*T_F
=I_{\mathcal H_n}.
\]
Thus $G_t$ is a dual frame of $F$
for every $t\in\mathbb R$.

Since $G$ is positive, there exists
$c=(c_i)_{i=1}^N\in\mathbb R_{++}^N$
such that
\[
T_G^*c=\sum_{i=1}^Nc_i g_i=0.
\]
Let $P_G$ denote the orthogonal
projection onto $\operatorname{Ran}T_G$.
For sufficiently small $|t|$, the
operators $T_{G_t}^*$ have full rank,
and their orthogonal kernel projections
$I-P_{G_t}$ converge to $I-P_G$.

Consequently,
\[
c_t=(I-P_{G_t})c\longrightarrow
(I-P_G)c=c.
\]
Since $c$ is strictly positive,
$c_t\in\mathbb R_{++}^N$ for all
sufficiently small $|t|$.
Moreover,
\[
T_{G_t}^*c_t=0.
\]
Therefore, $G_t$ is a positive
dual frame of $F$ whenever $|t|$
is sufficiently small.

Choose $t>0$ sufficiently small.
Since $H\neq0$, the positive dual
frames $G_t$ and $G_{-t}$ are distinct,
and
\[
G=\frac12G_t+\frac12G_{-t}.
\]
Hence $G$ is not an extreme point
of $\mathcal D_+(F)$.
Since $G$ was arbitrary, the proof
is complete.
\end{proof}

We next establish a geometric property of the
collection of positive dual frames. In particular,
we show that every positive dual frame of a
positive frame can be represented as a convex
combination of two relative boundary dual frames.

\begin{lemma}\label{lem:bounded-positive-duals}
Let $F=\{f_i\}_{i=1}^N$ be a positive frame
for $\mathcal H_n$, and let
$G=\{g_i\}_{i=1}^N\in\mathcal D_+(F)$.
Choose $d=(d_i)_{i=1}^N\in\ker T_F^*$
with $d_i>0$ for every $i$, and fix
$u\in\mathcal H_n\setminus\{0\}$.
Set $H=\{d_i u\}_{i=1}^N$ and
\[
A=\{s\in\mathbb R:
G+sH\in\mathcal D_+(F)\}.
\]
Then $A$ is a bounded open subset of
$\mathbb R$ containing $0$.
\end{lemma}

\begin{proof}
Since $d\in\ker T_F^*$, we have
\[
T_H^*T_Fx
=\sum_{i=1}^N\langle x,f_i\rangle d_i u
=\left\langle x,\sum_{i=1}^Nd_if_i
\right\rangle u=0
\]
for every $x\in\mathcal H_n$.
Consequently,
\[
T_{G+sH}^*T_F
=T_G^*T_F+sT_H^*T_F
=I_{\mathcal H_n},
\]
so $G+sH$ is a dual frame of $F$
for every $s\in\mathbb R$.

To establish boundedness, let
$x=u/\|u\|$ and define
\[
\alpha_i=
\frac{\langle x,g_i\rangle}{d_i\|u\|},
\qquad
m=\min_{1\leq i\leq N}\alpha_i,
\qquad
M=\max_{1\leq i\leq N}\alpha_i.
\]
For every $i$ and $s\in\mathbb R$,
\[
\langle x,g_i+sd_i u\rangle
=d_i\|u\|(\alpha_i+s).
\]

If $s\leq-M$, then
$\langle x,g_i+sd_i u\rangle\leq0$
for every $i$. By the positive-spanning
criterion, $G+sH$ is not positive.

Similarly, if $s\geq-m$, then
$\langle -x,g_i+sd_i u\rangle\leq0$
for every $i$, and hence $G+sH$
is not positive. Therefore,
\[
A\subseteq(-M,-m),
\]
which proves that $A$ is bounded.

Since positive spanning is an open
condition on finite families of vectors,
$A$ is open. Finally, $0\in A$ because
$G\in\mathcal D_+(F)$.
\end{proof}

The preceding lemma allows us to represent
every positive dual frame as a convex
combination of two relative boundary dual frames.

\begin{theorem}\label{thm:convex-boundary-duals}
Let $F=\{f_i\}_{i=1}^N$ be a positive frame
for $\mathcal H_n$. Then every positive
dual frame of $F$ can be expressed as a
convex combination of two distinct dual
frames belonging to the relative boundary
of $\mathcal D_+(F)$ in $\mathcal D(F)$.
\end{theorem}

\begin{proof}
Let $G\in\mathcal D_+(F)$. Since $F$
is positive, there exists
$d=(d_i)_{i=1}^N\in\ker T_F^*$
with $d_i>0$ for every $i$.
Fix $u\in\mathcal H_n\setminus\{0\}$
and set $H=\{d_i u\}_{i=1}^N$.

Consider
\[
A=\{s\in\mathbb R:
G+sH\in\mathcal D_+(F)\}.
\]
By Lemma~\ref{lem:bounded-positive-duals},
$A$ is bounded, open, and contains $0$.

Let $(a,b)$ be the connected component
of $A$ containing $0$. Since $A$ is
bounded and open, we have
\[
-\infty<a<0<b<\infty.
\]
Moreover, $G+sH\in\mathcal D_+(F)$
for every $s\in(a,b)$.

Define
\[
G_1=G+aH,\qquad G_2=G+bH.
\]
Since $G+sH$ is a dual frame of $F$
for every $s\in\mathbb R$, both
$G_1$ and $G_2$ belong to
$\mathcal D(F)$.

By the maximality of $(a,b)$ and
the openness of $A$, neither $a$
nor $b$ belongs to $A$.
Thus $G_1,G_2\notin\mathcal D_+(F)$.
On the other hand,
\[
G+sH\longrightarrow G_1
\quad\text{as }s\to a^+,
\]
and
\[
G+sH\longrightarrow G_2
\quad\text{as }s\to b^-.
\]
Consequently,
\[
G_1,G_2\in
\partial_{\mathcal D(F)}\mathcal D_+(F).
\]

Finally, since $a<0<b$, we have
\[
G=
\frac{b}{b-a}G_1
+\frac{-a}{b-a}G_2.
\]
Both coefficients are strictly positive
and their sum is one. Moreover,
$G_1\neq G_2$, since $a\neq b$
and $H\neq0$. Therefore, $G$ is a strict convex
combination of two distinct relative
boundary dual frames.
\end{proof}

\begin{remark}
The converse of Theorem~\ref{thm:convex-boundary-duals}
does not hold in general. In particular, a
convex combination of two relative boundary
dual frames need not be positive.

To see this, consider the positive frame
$F=\{e_1,-e_1,e_2,-e_2\}$ for $\mathbb R^2$,
where $\{e_1,e_2\}$ is the standard
orthonormal basis. Let
\[
G_1=\{e_1,0,e_2,0\},
\qquad
G_2=\{0,-e_1,e_2,0\}.
\]
Clearly, both $G_1$ and $G_2$ are dual
frames of $F$, but neither is positive.
Moreover, both belong to the relative
boundary of $\mathcal D_+(F)$ in
$\mathcal D(F)$, as they can be approximated
by positive dual frames.

For every $t\in(0,1)$, we have
\[
tG_1+(1-t)G_2
=
\{te_1,-(1-t)e_1,e_2,0\}.
\]
Since all vectors in this family have
nonnegative second coordinates, the family
does not positively span $\mathbb R^2$.
Consequently,
\[
tG_1+(1-t)G_2\notin\mathcal D_+(F),
\qquad 0<t<1.
\]
Thus, although every positive dual frame
of a positive frame is a convex combination
of two suitably chosen relative boundary
dual frames, the same conclusion does not
hold for arbitrary pairs of boundary duals.
\end{remark}

\begin{example}
Let $\{e_1,e_2\}$ be the standard orthonormal
basis of $\mathbb R^2$, and consider the frame
\[
F=\{e_1,e_2,e_1\}.
\]
Clearly, $F$ is not a positive frame.  Every dual frame of $F$ can be written as
\[
G(u)=\{e_1-u,e_2,u\},
\qquad u=(a,b)\in\mathbb R^2.
\]
Indeed, the duality condition is equivalent to
$g_1+g_3=e_1$ and $g_2=e_2$.

We claim that $G(u)$ is a positive dual
if and only if
\[
(a<0,\ b<0)
\quad\text{or}\quad
(a>1,\ b>0).
\]
To verify this, observe that a strictly
positive relation among the vectors of
$G(u)$ is equivalent to
\[
\alpha(e_1-u)+\beta e_2+\gamma u=0,
\qquad \alpha,\beta,\gamma>0.
\]
Comparing coordinates yields
\[
\alpha+(\gamma-\alpha)a=0,
\qquad
\beta+(\gamma-\alpha)b=0.
\]
These equations admit strictly positive
solutions precisely in the two regions
described above.

Now suppose that $G(u)$ is positive.

If $a<0$ and $b<0$, set
\[
u_1=(2a,0),\qquad u_2=(0,2b).
\]
Then
\[
G(u)=\frac12G(u_1)+\frac12G(u_2).
\]
Both $G(u_1)$ and $G(u_2)$ belong to the
relative boundary of $\mathcal D_+(F)$
in $\mathcal D(F)$.

Similarly, if $a>1$ and $b>0$, set
\[
u_1=(1,2b),\qquad u_2=(2a-1,0).
\]
Again,
\[
G(u)=\frac12G(u_1)+\frac12G(u_2),
\]
where $G(u_1)$ and $G(u_2)$ are distinct
boundary dual frames. Thus every positive dual frame of $F$
is a convex combination of two distinct
dual frames belonging to the relative
boundary of $\mathcal D_+(F)$, even though
$F$ itself is not positive.
\end{example}

The boundary behavior of positive dual frames
raises a natural question concerning the
preservation of positivity under convex
combinations. Although convex combinations
of dual frames preserve the reconstruction
identity, they need not preserve positive
spanning. The following example illustrates
this phenomenon by showing that the midpoint
of two boundary points of the set of positive
duals may fail to be positive.

\begin{example}
\label{ex:boundary-positive-duals}
Consider the frame
$F=\{(1,0),(1,1),(0,1)\}$ for
$\mathbb R^2$. For $0\leq t<1$, define
$G_t=\{g_1^{(t)},g_2^{(t)},g_3^{(t)}\}$,
where
\[
\begin{aligned}
g_1^{(t)}&=\frac{1}{2-t}(1-t,-1),\\
g_2^{(t)}&=\frac{1}{2-t}(1,1),\\
g_3^{(t)}&=\frac{1}{2-t}(-1,1-t).
\end{aligned}
\]
It is straightforward to verify that
$G_t$ is a dual frame of $F$ and
\[
g_1^{(t)}+t g_2^{(t)}
+g_3^{(t)}=0.
\]
Thus $G_t$ is positive for every
$0<t<1$. However, $G_t$ converges
as $t\to0^+$ to the semi-positive
dual frame
\[
G_0=
\left\{
\left(\frac12,-\frac12\right),
\left(\frac12,\frac12\right),
\left(-\frac12,\frac12\right)
\right\},
\]
which is not positive since
$\ker T_{G_0}^*
=\operatorname{span}\{(1,0,1)\}$.

Now consider another dual frame
\[
H=\{(0,0),(1,0),(-1,1)\}.
\]
Since
$\ker T_H^*
=\operatorname{span}\{(1,0,0)\}$,
the frame $H$ is also semi-positive
but not positive. Both $G_0$ and
$H$ belong to the relative boundary
of $\mathcal D_+(F)$.

Their midpoint $M=\frac12(G_0+H)$
satisfies
\[
\begin{aligned}
M=\bigg\{&
\left(\frac14,-\frac14\right),
\left(\frac34,\frac14\right),
\left(-\frac34,\frac34\right)
\bigg\}.
\end{aligned}
\]
A direct computation gives
\[
\ker T_M^*
=\operatorname{span}\{(3,0,1)\}.
\]
Hence $M$ is not positive.
This shows that positive dual frames
may converge to non-positive boundary
points and that convex combinations
of such boundary points need not
be positive.
\end{example}

\section{Positive Frames Arising from Finite Group Actions}
\label{sec:group-positive-frames-erasures}

We now investigate positive frames arising from finite
group representations, with particular emphasis on the
relationship between linear connectedness and
balancedness. The symmetry induced by a group action
provides additional information about the linear
relations among the frame vectors. In particular, the
relation space of a group-generated frame is invariant
under the natural permutation action of the group.
This observation allows us to establish a connection
between linear connectedness and the balancedness
of frames with excess one.

\subsection{\bf{Group Averaging, Balancedness, and Positivity}}
We investigate the relationship between balancedness
and positivity for frames generated by finite group
representations. Let $G$ be a finite group and let
$\pi:G\to O(\mathbb R^n)$ be an orthogonal representation.
Suppose that $\pi$ admits a decomposition into pairwise
inequivalent irreducible representations of the form
\[
\pi\simeq
\pi_1^{(m_1)}\oplus\cdots\oplus\pi_k^{(m_k)},
\qquad
\pi_i\not\simeq\pi_j \quad (i\neq j).
\]
For each $\xi\in\mathbb R^n$, consider the
group-generated family
$F_\xi=\{\pi(g)\xi:g\in G\}$, indexed by the
elements of $G$. Recall that $F_\xi$ is balanced
if $\sum_{g\in G}\pi(g)\xi=0$. The averaging
operator associated with $\pi$ is closely related
to the subspace of vectors fixed by the group
action. In particular, the absence of the trivial
representation determines when every orbit family
is balanced. Moreover, whenever such an orbit
family spans $\mathbb R^n$, its balancedness
implies positivity. The following proposition
establishes these equivalences.

\begin{proposition}
\label{prop:group-representation-balanced-positive}
Under the preceding assumptions, the following
statements are equivalent:
\begin{enumerate}
\item[(i)] The trivial representation does not
occur among $\pi_1,\ldots,\pi_k$.
\item[(ii)] $\displaystyle\sum_{g\in G}\pi(g)=0$.
\item[(iii)] $F_\xi$ is balanced for every
$\xi\in\mathbb R^n$.
\item[(iv)] $F_\xi$ is a positive frame whenever
it spans $\mathbb R^n$.
\end{enumerate}
\end{proposition}
\begin{proof}
Define the averaging operator
\[
P=\frac{1}{|G|}\sum_{g\in G}\pi(g).
\]
Let
\[
\mathcal I(\pi)
=\{x\in\mathbb R^n:\pi(g)x=x
\text{ for all }g\in G\}
\]
denote the subspace of invariant vectors.

For every $h\in G$, we have
\[
\pi(h)P
=\frac{1}{|G|}\sum_{g\in G}\pi(hg)
=P.
\]
Consequently, $\operatorname{ran}P
\subseteq\mathcal I(\pi)$. Conversely,
if $x\in\mathcal I(\pi)$, then
\[
Px=\frac{1}{|G|}
\sum_{g\in G}\pi(g)x=x.
\]
Thus $\mathcal I(\pi)\subseteq
\operatorname{ran}P$, and hence
\[
\operatorname{ran}P=\mathcal I(\pi).
\]
Since $\pi$ is orthogonal, $P$ is the
orthogonal projection onto $\mathcal I(\pi)$.

We first prove (i) $\Leftrightarrow$ (ii).
The trivial representation occurs in
$\pi$ if and only if
$\mathcal I(\pi)\neq\{0\}$. Therefore,
\[
\begin{aligned}
\text{(i)}
&\Longleftrightarrow
\mathcal I(\pi)=\{0\}\\
&\Longleftrightarrow P=0\\
&\Longleftrightarrow
\sum_{g\in G}\pi(g)=0,
\end{aligned}
\]
which proves the equivalence.

Next, (ii) $\Leftrightarrow$ (iii)
follows immediately from
\[
\sum_{g\in G}\pi(g)\xi
=\left(\sum_{g\in G}\pi(g)\right)\xi.
\]
Indeed, the right-hand side vanishes for
every $\xi\in\mathbb R^n$ if and only if
$\sum_{g\in G}\pi(g)=0$.

We now prove (iii) $\Rightarrow$ (iv).
Let $\xi\in\mathbb R^n$ be such that
$F_\xi$ spans $\mathbb R^n$.
By (iii),
\[
\sum_{g\in G}\pi(g)\xi=0.
\]
Thus $F_\xi$ admits a linear dependence
with all coefficients strictly positive.
Since $F_\xi$ spans $\mathbb R^n$,
it follows that $F_\xi$ positively spans
$\mathbb R^n$. Hence $F_\xi$ is a
positive frame.

Finally, suppose (iv) holds. We show
that $\mathcal I(\pi)=\{0\}$.
Assume, to the contrary, that
$\mathcal I(\pi)\neq\{0\}$. Then $P\neq0$.

If there exists a frame vector $\xi$,
then $F_\xi$ is positive by (iv).
Consequently, there exist coefficients
$c_g>0$, $g\in G$, such that
\[
\sum_{g\in G}c_g\pi(g)\xi=0.
\]
Applying $P$ and using
$P\pi(g)=P$ for every $g\in G$, we obtain
\[
0=\sum_{g\in G}c_gP\pi(g)\xi
=\left(\sum_{g\in G}c_g\right)P\xi.
\]
Since $\sum_{g\in G}c_g>0$, it follows
that $P\xi=0$. However, $P$ commutes
with every $\pi(g)$, so
\[
P\pi(g)\xi=\pi(g)P\xi=0
\qquad(g\in G).
\]
Since $F_\xi$ spans $\mathbb R^n$,
we conclude that $P=0$, a contradiction.

Therefore $\mathcal I(\pi)=\{0\}$,
which proves (i) and completes the proof.
\end{proof}

\subsection{\bf{Linear Connectedness and Excess-One Orbit Frames}}

We now investigate the relationship between
linear connectedness and balancedness for
frames generated by finite group representations.
The notion of linear connectedness was studied
by Pehlivan, Han, and Mohapatra
\cite{PehlivanHanMohapatra2013} in connection
with spectrally optimal dual frames for erasures.
We first recall the relevant definitions and
a characterization from their work.

Let $F=\{f_i\}_{i=1}^N$ be a finite family
of nonzero vectors in $\mathcal H_n$.
Two vectors $f_i$ and $f_j$ are said to be
\emph{linearly $F$-connected} if there exist
vectors $f_{k_1},\ldots,f_{k_m}$ from $F$
such that $\{f_j,f_{k_1},\ldots,f_{k_m}\}$
is linearly independent and
\[
f_i=cf_j+\sum_{\ell=1}^m c_\ell f_{k_\ell},
\]
where $c,c_1,\ldots,c_m$ are all nonzero.
The family $F$ is called \emph{linearly
connected} if every two of its vectors
are linearly $F$-connected.

For $\Lambda\subseteq\{1,\ldots,N\}$,
write
$\mathcal H_\Lambda
=\operatorname{span}\{f_i:i\in\Lambda\}$.
The family $F$ is said to have the
\emph{intersection dependent property}
if
\[
\mathcal H_\Lambda\cap
\mathcal H_{\Lambda^c}\neq\{0\}
\]
for every nonempty proper subset
$\Lambda\subset\{1,\ldots,N\}$.
Furthermore, $F$ is called
\emph{$k$-independent} if every
$k$-element subfamily is linearly
independent.

The following characterization was
established in
\cite[Theorem~2.1]{PehlivanHanMohapatra2013}.

\begin{theorem}\cite{PehlivanHanMohapatra2013}
\label{thm:linear-connectedness-characterization}
Let $F=\{f_i\}_{i=1}^N$ be a finite
family of nonzero vectors in
$\mathcal H_n$, and let
$d=\dim\operatorname{span}(F)$.
Then the following are equivalent:
\begin{enumerate}
\item[(i)] $F$ is linearly connected.
\item[(ii)] $F$ has the intersection
dependent property.
\item[(iii)] $F$ contains a
$d$-independent subfamily consisting
of at least $d+1$ vectors.
\end{enumerate}
\end{theorem}

We now consider the special case of
group-generated frames with excess one.
In this setting, linear connectedness
imposes a particularly strong restriction
on the unique linear dependence among
the orbit vectors. The following
proposition shows that this dependence
must be balanced.

\begin{proposition}
\label{prop:linearly-connected-balanced}
Let $G$ be a finite group and let
$\pi:G\to GL(\mathcal H_n)$ be a representation
of $G$ on an $n$-dimensional Hilbert space
$\mathcal H_n$. Suppose that
$F=\{\pi(g)\xi:g\in G\}$ is a linearly
connected frame for $\mathcal H_n$ with
$|G|=n+1$. Then $F$ is balanced, that is,
\[
\sum_{g\in G}\pi(g)\xi=0.
\]
\end{proposition}

\begin{proof}

Since $F$ is a frame for $\mathcal H_n$,
we have $\operatorname{rank}(T_F^*)=n.$
As $|G|=n+1$, the rank-nullity theorem gives $\dim\ker(T_F^*)=1.$ 

For each $h\in G$, define
$L_h:\mathbb C^{|G|}\to\mathbb C^{|G|}$ by
\[
(L_hc)_g=c_{h^{-1}g},
\qquad g\in G.
\]
Then, for every $c=(c_g)_{g\in G}$,
\[
\begin{aligned}
T_F^*(L_hc)
&=\sum_{g\in G}(L_hc)_g\pi(g)\xi\\
&=\sum_{g\in G}c_{h^{-1}g}\pi(g)\xi\\
&=\sum_{r\in G}c_r\pi(hr)\xi\\
&=\pi(h)\sum_{r\in G}c_r\pi(r)\xi\\
&=\pi(h)T_F^*c.
\end{aligned}
\]
Consequently,
\[
T_F^*L_h=\pi(h)T_F^*.
\]

If $c\in\ker(T_F^*)$, then
\[
T_F^*(L_hc)
=\pi(h)T_F^*c=0.
\]
Thus,
\[
L_h( kerT_F^*)\subseteq
 kerT_F^*.
\]
Since $L_h$ is invertible, we obtain
\[
L_h( kerT_F^*)= kerT_F^*.
\]

Now let $0\neq c=(c_g)_{g\in G}
\in\ker(T_F^*)$. Since the kernel
is one-dimensional, for every $h\in G$
there exists a nonzero scalar $\alpha_h$
such that
\[
L_hc=\alpha_hc.
\]
Equivalently,
\[
c_{h^{-1}g}=\alpha_hc_g,
\qquad g,h\in G.
\]

Since $F$ is linearly connected and
$\dim\ker(T_F^*)=1$, the coordinates
of $c$ are equal. Hence there exists
$a\neq0$ such that
\[
c=a\mathbf 1,
\]
where $\mathbf 1=(1,\ldots,1)^T$.
Therefore,
\[
0=T_F^*c
=a\sum_{g\in G}\pi(g)\xi.
\]
Since $a\neq0$, it follows that
\[
\sum_{g\in G}\pi(g)\xi=0.
\]
Hence $F$ is balanced.
\end{proof}

The group action imposes additional restrictions
on the linear relations among orbit vectors.
In the excess-one case, these relations are
completely determined by a real character
of the underlying group.

\begin{proposition}
\label{prop:orbit-kernel-character}
Let $G$ be a finite group and let
$\pi:G\to GL(\mathcal H_n)$ be a
representation on a real Hilbert space.
Suppose that
$F=\{\pi(g)\xi:g\in G\}$ is a frame
for $\mathcal H_n$ with $|G|=n+1$.
Then there exists a homomorphism
$\chi:G\to\{-1,1\}$ such that
\[
\ker T_F^*
=\operatorname{span}
\{(\chi(g))_{g\in G}\}.
\]
In particular, $F$ is linearly connected.
\end{proposition}

\begin{proof}
Since $|G|=n+1$, the synthesis kernel
$K=\ker T_F^*$ is one-dimensional.
For each $h\in G$, let $L_h$ denote
the permutation operator on
$\mathbb R^G$ defined by
$L_he_g=e_{hg}$.
Then
\[
T_F^*L_h=\pi(h)T_F^*,
\]
and consequently $L_hK=K$.

Choose a nonzero $c\in K$. Since
$K$ is one-dimensional, there exists
a scalar $\chi(h)\neq0$ such that
$L_hc=\chi(h)c$.
The identity $L_{hk}=L_hL_k$
implies that $\chi$ is a homomorphism.
Since $G$ is finite and $\chi$ is
real-valued, $\chi(h)\in\{-1,1\}$
for every $h\in G$.

Moreover, the transitivity of the
permutation action implies that no
coordinate of $c$ can vanish.
Using the coordinate corresponding
to the identity element, we obtain
$c_g=c_e\chi(g^{-1})=c_e\chi(g)$.
Thus
\[
K=\operatorname{span}
\{(\chi(g))_{g\in G}\}.
\]
Since the unique linear dependence
has no vanishing coefficient, every
$n$-element subfamily of $F$ is
linearly independent. Hence $F$
is linearly connected by
\cite[Theorem~2.1]{PehlivanHanMohapatra2013}.
\end{proof}

The preceding proposition yields a
characterization of balancedness and
positive spanning for excess-one orbit
frames.

\begin{corollary}
\label{thm:orbit-balanced-positive}
Under the hypotheses of
Proposition~\ref{prop:orbit-kernel-character},
the following statements are equivalent:
\begin{enumerate}
\item[(i)] $F$ is balanced.
\item[(ii)] $F$ is a positive frame.
\item[(iii)] The character $\chi$
is trivial.
\item[(iv)] $\mathcal H_n$ contains
no nonzero vector fixed by every
$\pi(g)$, $g\in G$.
\end{enumerate}
\end{corollary}

\begin{proof}
By Proposition~\ref{prop:orbit-kernel-character},
the synthesis kernel is generated by
$(\chi(g))_{g\in G}$.
Therefore, the all-ones vector belongs
to $\ker T_F^*$ if and only if
$\chi$ is trivial. This proves
(i)$\Leftrightarrow$(iii).

Since $\chi(g)\in\{-1,1\}$, the
kernel contains a strictly positive
vector if and only if $\chi$ is
trivial. By Proposition~3.1,
(ii)$\Leftrightarrow$(iii).

Finally, the averaging operator
\[
P=\frac1{|G|}\sum_{g\in G}\pi(g)
\]
is a projection onto the common
fixed-point subspace of the
representation. Since $\xi$ is
cyclic, the range of $P$ is generated
by $P\xi$. Consequently, the
fixed-point subspace is trivial
if and only if $P\xi=0$, which is
equivalent to balancedness.
\end{proof}

An immediate consequence is that
excess-one orbit frames generated
by groups of odd order are necessarily
positive.

\begin{corollary}
\label{cor:odd-order-positive-orbit}
Let $G$ be a finite group of odd
order and let
$F=\{\pi(g)\xi:g\in G\}$ be a frame
for $\mathcal H_n$ with $|G|=n+1$.
Then $F$ is balanced and positive.
In particular, $F$ is a minimal
positive basis for $\mathcal H_n$.
\end{corollary}

\begin{example}
\label{ex:connected-nonbalanced-orbit}
Let $G=C_4=\langle a:a^4=e\rangle$
and consider the representation
$\pi:G\to GL(\mathbb R^3)$ determined by
\[
\pi(a)=
\begin{pmatrix}
1&0&0\\
0&0&-1\\
0&1&0
\end{pmatrix}.
\]
For $\xi=(1,1,0)$, the corresponding
orbit frame is
\[
F=\{(1,1,0),(1,0,1),
(1,-1,0),(1,0,-1)\}.
\]
A direct computation gives
\[
\ker T_F^*
=\operatorname{span}\{(1,-1,1,-1)\}.
\]
Thus every three vectors of $F$
are linearly independent, and hence
$F$ is linearly connected. However,
\[
\sum_{g\in G}\pi(g)\xi=(4,0,0)\neq0.
\]
Therefore $F$ is neither balanced
nor positive.

The corresponding character is
$\chi(a^j)=(-1)^j$, illustrating
the obstruction to balancedness
described in
Theorem~\ref{thm:orbit-balanced-positive}.
\end{example}

\section{Positive Dependence, Convex Geometry, and Erasure Robustness}
\label{sec:positive-dependence-erasures}

We begin by recalling the notion of minimal positive
dependence in finite-dimensional real Hilbert spaces.
Let $I$ be a finite index set and let
$F_I=\{f_i:i\in I\}$ be a family of nonzero vectors
in $\mathcal H_n$, with
$S_I=\operatorname{span}\{f_i:i\in I\}$.
The family $F_I$ is said to be
\emph{positively dependent} if there exist
strictly positive scalars $c_i$, $i\in I$, such that
$\sum_{i\in I}c_i f_i=0$. It is called
\emph{minimally positively dependent} if it is
positively dependent, but no proper nonempty
subfamily admits a nontrivial nonnegative linear
dependence, that is, a relation with all coefficients
nonnegative and at least one coefficient strictly
positive.

\subsection{\bf{Convex Geometry of Positive Relations}}
Suppose that $F=\{f_i\}_{i=1}^N$
is an $(\ell-1)$-redundant and $\ell$-exact
positive frame for $\mathbb R^n$, where
$\ell\geq 1$. By the known lower bound for
redundant positive frames, we have
$N\geq n+2\ell-1$. On the other hand,
every positive frame contains a positive
basis $B$ with $|B|\leq 2n$. Since $F$
is $\ell$-exact, there can be at most
$\ell-1$ vectors outside $B$. Consequently,
\[
n+2\ell-1\leq N\leq 2n+\ell-1.
\]
In particular, $\ell\leq n$.

For $\ell=1$, both bounds are attained.
Indeed, the family
\[
F_1=\left\{e_1,\ldots,e_n,
-\sum_{i=1}^n e_i\right\}
\]
is an exact positive frame for $\mathbb R^n$
with $n+1$ elements, attaining the lower
bound. Similarly,
\[
F_2=\{\pm e_i:1\leq i\leq n\}
\]
is an exact positive frame with $2n$
elements, attaining the upper bound.
Here $\{e_i\}_{i=1}^n$ denotes the
standard orthonormal basis of $\mathbb R^n$.

Let $F=\{f_i\}_{i=1}^N$ be a finite family of vectors in
a real Hilbert space $\mathcal H_n$. Consider the convex set
\[
\Omega=\left\{
\lambda=(\lambda_1,\ldots,\lambda_N)\in[0,\infty)^N:
\sum_{i=1}^N\lambda_i=1,\quad
\sum_{i=1}^N\lambda_i f_i=0
\right\}.
\]
For each $\lambda\in\Omega$, define
$I(\lambda)=\{i:\lambda_i>0\}$. The following proposition characterizes the extreme points
of $\Omega$ in terms of linear independence of the
corresponding augmented vectors.

\begin{proposition}
Let $\lambda\in\Omega$. Then $\lambda$ is an extreme
point of $\Omega$ if and only if the family
\[
\{(f_i,1):i\in I(\lambda)\}
\]
is linearly independent in $\mathcal H_n\oplus\mathbb R$.
\end{proposition}

\begin{proof}
Suppose first that
$\{(f_i,1):i\in I(\lambda)\}$ is linearly independent.
Assume that
\[
\lambda=t\mu+(1-t)\nu,
\qquad 0<t<1,
\]
for some $\mu,\nu\in\Omega$.

If $j\notin I(\lambda)$, then $\lambda_j=0$. Since
$\mu_j,\nu_j\geq0$, we have
\[
0=t\mu_j+(1-t)\nu_j,
\]
which implies $\mu_j=\nu_j=0$.

Consequently, using $\mu,\nu\in\Omega$, we obtain
\[
\sum_{i\in I(\lambda)}(\mu_i-\nu_i)f_i=0
\quad\text{and}\quad
\sum_{i\in I(\lambda)}(\mu_i-\nu_i)=0.
\]
Equivalently,
\[
\sum_{i\in I(\lambda)}
(\mu_i-\nu_i)(f_i,1)=0.
\]
By linear independence,
$\mu_i=\nu_i$ for every $i\in I(\lambda)$.
Together with the equality outside $I(\lambda)$,
this gives $\mu=\nu=\lambda$.
Hence $\lambda$ is an extreme point of $\Omega$.

Conversely, suppose that $\lambda$ is an extreme
point of $\Omega$. Assume, to the contrary, that
$\{(f_i,1):i\in I(\lambda)\}$ is linearly dependent.
Then there exist scalars $a_i$, not all zero, such that
\[
\sum_{i\in I(\lambda)}a_i f_i=0
\quad\text{and}\quad
\sum_{i\in I(\lambda)}a_i=0.
\]
Extend $a=(a_1,\ldots,a_N)$ by setting $a_i=0$
for $i\notin I(\lambda)$.

Since $\lambda_i>0$ for every $i\in I(\lambda)$,
we can choose $\varepsilon>0$ sufficiently small
so that
\[
\lambda_i\pm\varepsilon a_i\geq0
\qquad (1\leq i\leq N).
\]
Define
\[
\lambda^+=\lambda+\varepsilon a,
\qquad
\lambda^-=\lambda-\varepsilon a.
\]
Then $\lambda^+,\lambda^-\in[0,\infty)^N$.
Moreover,
\[
\sum_{i=1}^N\lambda_i^\pm
=\sum_{i=1}^N\lambda_i
\pm\varepsilon\sum_{i=1}^N a_i=1
\]
and
\[
\sum_{i=1}^N\lambda_i^\pm f_i
=\sum_{i=1}^N\lambda_i f_i
\pm\varepsilon\sum_{i=1}^N a_i f_i=0.
\]
Thus $\lambda^+,\lambda^-\in\Omega$.

Since $a\neq0$, we have $\lambda^+\neq\lambda^-$,
while
\[
\lambda=\frac12\lambda^+
+\frac12\lambda^-.
\]
This contradicts the assumption that $\lambda$
is an extreme point of $\Omega$.
Therefore $\{(f_i,1):i\in I(\lambda)\}$ is
linearly independent.
\end{proof}

\begin{corollary}
Let $F=\{f_i\}_{i=1}^N$ be a finite family in an
$n$-dimensional real Hilbert space $\mathcal H$.
If $\lambda$ is an extreme point of $\Omega$, then
\[
|I(\lambda)|\leq n+1.
\]
\end{corollary}

\begin{proof}
By the preceding proposition, the family
$\{(f_i,1):i\in I(\lambda)\}$ is linearly
independent in $\mathcal H\oplus\mathbb R$.
Since $\dim(\mathcal H\oplus\mathbb R)=n+1$,
the conclusion follows.
\end{proof}

\begin{corollary}
Let $\lambda$ be an extreme point of $\Omega$.
Then $\lambda$ is the unique element of $\Omega$
whose support is contained in $I(\lambda)$.
\end{corollary}

\begin{proof}
Suppose $\mu\in\Omega$ satisfies
$I(\mu)\subseteq I(\lambda)$. Since
$\lambda,\mu\in\Omega$, we have
\[
\sum_{i\in I(\lambda)}(\lambda_i-\mu_i)f_i=0
\quad\text{and}\quad
\sum_{i\in I(\lambda)}(\lambda_i-\mu_i)=0.
\]
Consequently,
\[
\sum_{i\in I(\lambda)}
(\lambda_i-\mu_i)(f_i,1)=0.
\]
By the preceding proposition,
$\{(f_i,1):i\in I(\lambda)\}$ is linearly
independent. Hence $\lambda_i=\mu_i$
for every $i\in I(\lambda)$, and therefore
$\lambda=\mu$.
\end{proof}

Let $F=\{f_i\}_{i=1}^N$ be a finite frame
consisting of pairwise distinct vectors in a
real Hilbert space $\mathcal H_n$, and let
$T_F$ and $T_F^*$ denote its analysis
and synthesis operators, respectively. The following
result characterizes the extreme points of
$\operatorname{conv}(F)$ in terms of the kernel
of the synthesis operator.

\begin{proposition}
\label{prop:extreme-points-kernel}
For each $i\in\{1,\ldots,N\}$,
$f_i\in\operatorname{Ext}(\operatorname{conv}(F))$
if and only if there does not exist
$d=(d_1,\ldots,d_N)\in kerT_F^*$
such that
\[
d_i<0,\qquad d_j\geq0\quad(j\neq i),
\qquad \sum_{j=1}^N d_j=0.
\]
\end{proposition}

\begin{proof}
Suppose that
$f_i\notin\operatorname{Ext}(\operatorname{conv}(F))$.
Since the vectors of $F$ are pairwise distinct,
there exist scalars $\lambda_j\geq0$, $j\neq i$,
such that
\[
f_i=\sum_{j\neq i}\lambda_jf_j,
\qquad \sum_{j\neq i}\lambda_j=1.
\]
Define $d=(d_1,\ldots,d_N)$ by
$d_i=-1$ and $d_j=\lambda_j$ for $j\neq i$.
Then
\[
T_F^*d
=-f_i+\sum_{j\neq i}\lambda_jf_j=0.
\]
Thus $d\in kerT_F^*$, with
$d_i<0$, $d_j\geq0$ for every $j\neq i$, and
$\sum_{j=1}^N d_j=0$.

Conversely, suppose that there exists
$d\in kerT_F^*$ satisfying
\[
d_i<0,\qquad d_j\geq0\quad(j\neq i),
\qquad \sum_{j=1}^N d_j=0.
\]
Since $T_F^*d=0$, we have
\[
d_if_i+\sum_{j\neq i}d_jf_j=0.
\]
Consequently,
\[
f_i=\sum_{j\neq i}\frac{d_j}{-d_i}f_j,
\qquad
\sum_{j\neq i}\frac{d_j}{-d_i}=1.
\]
Since all coefficients are nonnegative,
$f_i\in\operatorname{conv}\{f_j:j\neq i\}$.
As the frame vectors are pairwise distinct,
$f_i$ cannot be an extreme point of
$\operatorname{conv}(F)$.
\end{proof}

\subsection{\bf{Positive Spanning and Erasure Robustness}}
We next introduce a quantity that measures the
robustness of a positive frame under arbitrary
vector erasures. Let $F=\{f_i\}_{i=1}^N$ be a
finite family of vectors in $\mathbb R^n$.
For each $0\neq x\in\mathbb R^n$, define
\[
N_F^+(x)
=\#\{i\in\{1,\ldots,N\}:
\langle x,f_i\rangle>0\},
\]
where $\#$ denotes the cardinality of a set.
We further define
\[
d_+(F)=\min_{x\in\mathbb R^n\setminus\{0\}}
N_F^+(x).
\]
Equivalently, by homogeneity,
\[
d_+(F)=\min_{\|x\|=1}N_F^+(x).
\]

Recall that $F$ is a positive frame for
$\mathbb R^n$ if
$\operatorname{pos}(F)=\mathbb R^n$, where
\[
\operatorname{pos}(F)
=\left\{\sum_{i=1}^N a_i f_i:a_i\geq0\right\}.
\]
By the positive spanning criterion, this is
equivalent to the requirement that, for every
$0\neq x\in\mathbb R^n$, there exists an
index $i$ such that $\langle x,f_i\rangle>0$.

For a nonnegative integer $s$, we say that
$F$ is \emph{$s$-redundant as a positive frame}
if the deletion of any collection of at most
$s$ vectors leaves a positive frame for
$\mathbb R^n$.

\begin{proposition}
\label{prop:positive-erasure-characterization}
Let $F=\{f_i\}_{i=1}^N$ be a finite family
of vectors in $\mathbb R^n$, and let
$s\in\{0,1,\ldots,N\}$. Then $F$ remains
a positive frame after the deletion of any
collection of at most $s$ vectors if and
only if
\[
d_+(F)\geq s+1.
\]
\end{proposition}

\begin{proof}
Suppose first that $d_+(F)\geq s+1$.
Let $J\subseteq\{1,\ldots,N\}$ be arbitrary
with $|J|\leq s$, and consider the remaining
family
\[
F_{J^c}=\{f_i:i\notin J\}.
\]
For every $0\neq x\in\mathbb R^n$, there
are at least $s+1$ indices $i$ such that
$\langle x,f_i\rangle>0$. Since at most
$s$ indices have been removed, there exists
at least one index $i\notin J$ satisfying
\[
\langle x,f_i\rangle>0.
\]
By the positive spanning criterion,
$F_{J^c}$ positively spans $\mathbb R^n$.
As $J$ was arbitrary, $F$ remains positive
after any deletion of at most $s$ vectors.

Conversely, suppose that $F$ remains a
positive frame after the deletion of any
collection of at most $s$ vectors.
Assume, to the contrary, that
$d_+(F)\leq s$. Then there exists
$0\neq x\in\mathbb R^n$ such that
\[
\#\{i:\langle x,f_i\rangle>0\}\leq s.
\]
Let
\[
J=\{i:\langle x,f_i\rangle>0\}.
\]
Clearly, $|J|\leq s$. Upon deleting
the vectors indexed by $J$, we obtain
\[
\langle x,f_i\rangle\leq0
\qquad\text{for every }i\notin J.
\]
Thus the remaining family fails the
positive spanning criterion and hence
is not a positive frame. This contradicts
our assumption. Therefore $d_+(F)\geq s+1$, completing
the proof.
\end{proof}

\subsection{\bf{Exactness and Cardinality Bounds}}

We now investigate the relationship between positive
spanning, exactness, and the cardinality of positive
frames. In particular, we examine how the removal of
frame vectors affects the positive spanning property
and characterize positive frames that lose this
property upon the deletion of a prescribed number
of vectors. We also derive cardinality bounds that
reflect the interplay between redundancy and
exactness, providing further insight into the
combinatorial structure of positive frames.

\begin{proposition}\label{prop:s-redundant-positive-duals}
Let $F=\{f_i\}_{i=1}^N$ be a frame for
$\mathcal H_n$, and let $s$ be a nonnegative
integer with $N-s\geq n+1$.
Then, for $G=\{g_i\}_{i=1}^N\in\mathcal D(F)$,
the following statements are equivalent:
\begin{enumerate}
\item[(i)] $G$ is an $s$-redundant positive
frame.
\item[(ii)] For every nonzero
$x\in\mathcal H_n$,
\[
\left|\left\{
i:\langle x,g_i\rangle>0
\right\}\right|\geq s+1.
\]
\end{enumerate}
\end{proposition}

\begin{proof}
$(i)\Rightarrow(ii)$.
Suppose that $G$ is an $s$-redundant
positive frame. Assume, to the contrary,
that there exists a nonzero
$x\in\mathcal H_n$ such that
\[
\left|\left\{
i:\langle x,g_i\rangle>0
\right\}\right|\leq s.
\]
Set
\[
I_+(x)=\{i:\langle x,g_i\rangle>0\}.
\]
Choose a subset
$\Lambda\subseteq\{1,\ldots,N\}$
with $|\Lambda|=s$ and
$I_+(x)\subseteq\Lambda$.
Then
\[
\langle x,g_i\rangle\leq0
\qquad\text{for every }i\in\Lambda^c.
\]
By the positive-spanning criterion,
the subfamily $\{g_i\}_{i\in\Lambda^c}$
is not a positive frame, contradicting
the $s$-redundancy of $G$.
Hence (ii) holds.

$(ii)\Rightarrow(i)$.
Suppose that (ii) holds, and let
$\Lambda\subseteq\{1,\ldots,N\}$
be arbitrary with $|\Lambda|=s$.
For every nonzero $x\in\mathcal H_n$,
there are at least $s+1$ indices $i$
such that $\langle x,g_i\rangle>0$.
Since $|\Lambda|=s$, at least one
such index belongs to $\Lambda^c$.
Consequently,
\[
\max_{i\in\Lambda^c}
\langle x,g_i\rangle>0
\]
for every nonzero $x\in\mathcal H_n$.
By the positive-spanning criterion,
$\{g_i\}_{i\in\Lambda^c}$ positively
spans $\mathcal H_n$.

Since $\Lambda$ was arbitrary,
$G$ is an $s$-redundant positive frame.
\end{proof}

The preceding proposition also provides a
characterization of positive frames that are
$f$-redundant and $(f+1)$-exact. Here,
$(f+1)$-exactness means that deleting any
$f+1$ vectors destroys the positive spanning
property.

Combining the preceding characterization of exactness
with Proposition~\ref{prop:positive-erasure-characterization},
we obtain the following criterion for positive frames
that are simultaneously $f$-redundant and
$(f+1)$-exact.

\begin{proposition}
\label{thm:q-exact-characterization}
Let $F=\{f_i\}_{i=1}^N$ be a positive frame
for $\mathbb R^n$, and let $1\leq \ell \leq N$.
Then the following statements are equivalent:
\begin{enumerate}
\item[(i)] $F$ is $\ell$-exact.
\item[(ii)] For every $J\subseteq\{1,\ldots,N\}$
with $|J|=\ell$, there exists $x_J\neq 0$ such that
\[
\langle x_J,f_i\rangle\leq 0
\qquad (i\notin J).
\]
\item[(iii)] For every $J\subseteq\{1,\ldots,N\}$
with $|J|=\ell$, there exists $x_J\neq 0$ such that
\[
P_F(x_J)\subseteq J.
\]
\end{enumerate}
\end{proposition}

\begin{corollary}
\label{cor:f-redundant-exact}
Let $F=\{f_i\}_{i=1}^N$ be a positive frame
for $\mathbb R^n$, and let $0\leq f<N$.
Then $F$ is $f$-redundant and $(f+1)$-exact
if and only if $d_+(F)=f+1$
and, for every $J\subseteq\{1,\ldots,N\}$
with $|J|=f+1$, there exists $x_J\neq 0$
such that $P_F(x_J)=J.$
\end{corollary}

The preceding characterization yields an explicit
expression for the maximum number of erasures
under which a finite positive frame remains positive.

\begin{corollary}
\label{cor:max-positive-erasures}
Let $F=\{f_i\}_{i=1}^N$ be a finite positive frame
for a real Hilbert space $\mathcal H$. Then the
maximum number of arbitrary erasures that preserve
positivity is
\[
r_+(F)=d_+(F)-1.
\]
\end{corollary}

\begin{proof}
By Proposition~\ref{prop:positive-erasure-characterization},
$F$ remains positive after any deletion of at most
$s$ vectors precisely when $s<d_+(F)$.
Therefore, the largest such integer is
$d_+(F)-1$.
\end{proof}

\subsection{\bf{Approximation by Positive Parseval Frames}}

\noindent Let $\mathcal H_n$ be an $n$-dimensional real Hilbert space,
and let $N>n$ be fixed. Consider the set
\[
\mathcal K
=
\left\{
G=\{g_i\}_{i=1}^N:
G \text{ is a positive Parseval frame for }\mathcal H_n
\right\}.
\]
We equip $\mathcal H_{n}^N$ with the norm
$\|F\|_2=(\sum_{i=1}^N\|f_i\|^2)^{1/2}$.

We first observe that $\mathcal K$ is bounded but
neither convex nor closed. Indeed, since every
$G=\{g_i\}_{i=1}^N\in\mathcal K$ is Parseval, $\sum_{i=1}^N\|g_i\|^2
=
\operatorname{tr}(S_G)=n.$
Thus $\|G\|_2=\sqrt n$ for every $G\in\mathcal K$,
which shows that $\mathcal K$ is bounded.
On the other hand, if $G\in\mathcal K$, then
$-G\in\mathcal K$, whereas
$\frac12G+\frac12(-G)=0\notin\mathcal K$.
Hence $\mathcal K$ is not convex.

Furthermore, $\mathcal K$ is not closed.
To see this, let $\{e_1,\ldots,e_n\}$ be an
orthonormal basis of $\mathcal H_n$ and put
$N=n+1$. For $0<t<1$, define
\[
v_i(t)=e_i,\quad 1\leq i\leq n,
\qquad
v_{n+1}(t)=-t\sum_{i=1}^n e_i.
\]
Then $V_t=\{v_i(t)\}_{i=1}^{n+1}$ is a positive
frame for every $t>0$. Let $S_t$ denote its
frame operator and define
\[
G_t=\{S_t^{-1/2}v_i(t)\}_{i=1}^{n+1}.
\]
Since $S_t^{-1/2}$ is invertible, each $G_t$
is a positive Parseval frame. Moreover,
$S_t\to I_{\mathcal H_n}$ as $t\to0^+$, and hence
\[
G_t\longrightarrow
\{e_1,\ldots,e_n,0\},
\]
which is Parseval but not positive.
Thus $\mathcal K$ is not closed.
For $N>n+1$, the same conclusion follows
by adjoining $N-n-1$ zero vectors to each
frame in this construction.

We now determine the distance from an arbitrary
positive frame to the set $\mathcal K$.

\begin{proposition}\label{prop:positive-parseval-geometry}
For every positive frame
$F=\{f_i\}_{i=1}^N$ with frame operator $S_F$, we have
\[
\operatorname{dist}(F,\mathcal K)
=
\|S_F^{1/2}-I_{\mathcal H_n}\|_{\mathrm{HS}}
=
\left(
\sum_{j=1}^n(\sqrt{\lambda_j}-1)^2
\right)^{1/2},
\]
where $\lambda_1,\ldots,\lambda_n$ are the eigenvalues
of $S_F$. Furthermore, the distance is attained at
the canonical Parseval frame
$\{S_F^{-1/2}f_i\}_{i=1}^N$.
\end{proposition}

\begin{proof}
Let $A=\Theta_F^*$ be the synthesis operator
of $F$, so that $AA^*=S_F$. Define
$U=S_F^{-1/2}A$. Then
\[
UU^*
=
S_F^{-1/2}AA^*S_F^{-1/2}
=
I_{\mathcal H}.
\]
Since $F$ is positive and $S_F^{-1/2}$ is
invertible, the family
$\{S_F^{-1/2}f_i\}_{i=1}^N$ is a positive
Parseval frame. Thus $U$ is the synthesis
operator of an element of $\mathcal K$.

Let $G=\{g_i\}_{i=1}^N\in\mathcal K$ be arbitrary,
and let $B=\Theta_G^*$ be its synthesis operator.
Since $G$ is Parseval, $BB^*=I_{\mathcal H}$.
Using $A=S_F^{1/2}U$, we obtain
\[
\begin{aligned}
\|A-B\|_{\mathrm{HS}}^2
&=\operatorname{tr}(AA^*)
+\operatorname{tr}(BB^*)
-2\operatorname{tr}(AB^*)\\
&=\operatorname{tr}(S_F)+n
-2\operatorname{tr}(S_F^{1/2}UB^*).
\end{aligned}
\]
Since $UU^*=BB^*=I_{\mathcal H}$, we have
$\|U\|_{\mathrm{op}}=\|B^*\|_{\mathrm{op}}=1$.
Consequently,
\[
\|UB^*\|_{\mathrm{op}}
\leq
\|U\|_{\mathrm{op}}\|B^*\|_{\mathrm{op}}
=1.
\]
As $S_F^{1/2}$ is positive definite, it follows
that
\[
\operatorname{tr}(S_F^{1/2}UB^*)
\leq
\operatorname{tr}(S_F^{1/2}).
\]
Therefore,
\[
\begin{aligned}
\|A-B\|_{\mathrm{HS}}^2
\geq
\operatorname{tr}(S_F)+n
-2\operatorname{tr}(S_F^{1/2})
=
\|S_F^{1/2}-I_{\mathcal H_n}\|_{\mathrm{HS}}^2.
\end{aligned}
\]
Equality holds when $B=U$, corresponding
to the canonical Parseval frame
$\{S_F^{-1/2}f_i\}_{i=1}^N$. Hence
\[
\operatorname{dist}(F,\mathcal K)
=
\|S_F^{1/2}-I_{\mathcal H_n}\|_{\mathrm{HS}}.
\]
Finally, if $\lambda_1,\ldots,\lambda_n$
are the eigenvalues of $S_F$, then
\[
\|S_F^{1/2}-I_{\mathcal H_n}\|_{\mathrm{HS}}^2
=
\sum_{j=1}^n(\sqrt{\lambda_j}-1)^2,
\]
which completes the proof.
\end{proof}

We now consider the distance from an arbitrary frame, not
necessarily positive, to the set $\mathcal K$. Let
$F=\{f_i\}_{i=1}^N$ be a frame for $\mathcal H_n$,
with analysis operator $T_F$ and frame operator
$S_F=T_F^*T_F$. Recall that a frame
$G=\{g_i\}_{i=1}^N$ is positive if and only if
there exists a vector $c\in\mathbb R_{>0}^N$
such that $T_G^*c=0$.

For each $c\in\mathbb R^N$ with $\|c\|=1$,
let $P_c=I_N-cc^*$ denote the orthogonal
projection onto $c^\perp$. Since every positive
Parseval frame $G$ satisfies $T_G^*T_G=I_{\mathcal H_n}$
and $T_G^*c=0$ for some strictly positive unit
vector $c$, we have
\[
\begin{aligned}
\|T_F^*-T_G^*\|_{\mathrm{HS}}^2
=\operatorname{tr}(S_F)+n
-2\operatorname{tr}(T_F^*T_G)=\operatorname{tr}(S_F)+n
-2\operatorname{tr}(T_F^*P_cT_G).
\end{aligned}
\]
For a fixed $c>0$ with $\|c\|=1$, the
rectangular orthogonal Procrustes formula gives
\[
\sup_{\substack{T_G^*T_G=I_{\mathcal H_n}\\
T_G^*c=0}}
\operatorname{tr}(T_F^*P_cT_G)
=
\operatorname{tr}
\left((T_F^*P_cT_F)^{1/2}\right).
\]
Consequently, taking the infimum over all
positive Parseval frames and using the continuity
of the trace of the positive square root, we obtain
\[
\begin{aligned}
\operatorname{dist}(F,\mathcal K)^2
=\operatorname{tr}(S_F)+n-2\max_{\substack{c\in\mathbb R_{\geq0}^N
\|c\|=1}}
\operatorname{tr}
\left((T_F^*P_cT_F)^{1/2}\right).
\end{aligned}
\]
Since
$T_F^*P_cT_F=S_F-(T_F^*c)\otimes(T_F^*c)$,
this can equivalently be written as
\[\begin{aligned}
\operatorname{dist}(F,\mathcal K)^2
=\operatorname{tr}(S_F)+n-2\max_{\substack{c\in\mathbb R_{\geq0}^N
\|c\|=1}}
\operatorname{tr}\left(
\left[S_F-(T_F^*c)\otimes(T_F^*c)\right]^{1/2}
\right).
\end{aligned}
\]

Furthermore, since $0\leq P_c\leq I_N$,
we have $T_F^*P_cT_F\leq S_F$. By the
operator monotonicity of the square root,
\[
\operatorname{tr}
\left((T_F^*P_cT_F)^{1/2}\right)
\leq\operatorname{tr}(S_F^{1/2}).
\]
It follows that
\[
\operatorname{dist}(F,\mathcal K)
\geq
\|S_F^{1/2}-I_{\mathcal H_n}\|_{\mathrm{HS}}.
\]
Moreover, equality holds if and only if
there exists a nonzero vector
$c\in\mathbb R_{\geq0}^N$ such that
$T_F^*c=0$. In particular, the distance
formula in Proposition
\ref{prop:positive-parseval-geometry}
remains valid for every semi-positive frame
admitting a nontrivial nonnegative dependence.
On the other hand, if
$\ker T_F^*\cap\mathbb R_{\geq0}^N=\{0\}$,
then
\[
\operatorname{dist}(F,\mathcal K)
>
\|S_F^{1/2}-I_{\mathcal H_n}\|_{\mathrm{HS}}.
\]

To illustrate this distinction, consider the
frame $F=\{1/\sqrt2,1/\sqrt2\}$ for
$\mathcal H_1=\mathbb R$. Clearly,
$S_F=I_{\mathcal H_1}$, so $F$ is Parseval.
However, $F$ is not positive, since both
frame vectors are strictly positive. Every
positive Parseval frame with two vectors
has one positive and one negative entry.
Up to the choice of signs, such a frame can
be written as
\[
G_\theta=\{\cos\theta,-\sin\theta\},
\qquad 0<\theta<\frac{\pi}{2}.
\]
Therefore,
\[
\begin{aligned}
\|F-G_\theta\|_2^2
&=\left(\frac1{\sqrt2}-\cos\theta\right)^2
+\left(\frac1{\sqrt2}+\sin\theta\right)^2\\
&=2-\sqrt2(\cos\theta-\sin\theta).
\end{aligned}
\]
Taking the infimum over $0<\theta<\pi/2$
gives
\[
\operatorname{dist}(F,\mathcal K)
=\sqrt{2-\sqrt2}>0.
\]
The infimum is approached as $\theta\to0^+$,
but it is not attained, since the limiting
Parseval frame $\{1,0\}$ is not positive.
Thus, even when a frame is Parseval, its
distance to the set of positive Parseval
frames need not vanish.

\section{Concluding Remarks}

In this paper, we have investigated the geometric and
operator-theoretic structure of positive frames and
their duals in finite-dimensional real Hilbert spaces.
Our approach emphasizes the role of positive linear
dependence and the synthesis kernel in understanding
positive spanning, frame equivalence, and redundancy.
In particular, we have characterized equivalence
classes of positive frames through their relation
spaces, obtaining parametrizations in terms of the
relative interior of a standard simplex and suitable
subsets of Grassmannians.

A central part of our investigation concerns the
existence and structure of positive dual frames.
We have established that a finite frame admits a
positive dual precisely when it is redundant and
provided an explicit construction using a rank-one
perturbation of the canonical dual. We have also
examined the geometry of the collection of positive
duals, including its relative openness and the absence
of extreme points, and obtained criteria for the
existence of dual frames satisfying prescribed
positive relations. Furthermore, we have considered
positive-basis duals and the boundary behavior of
semi-positive frames.

We have also studied positive frames generated by
finite group representations, highlighting the
connections between balancedness, linear
connectedness, and positive spanning. In addition,
our investigation of normalized nonnegative relation
spaces provides a convex-geometric description of
positive dependencies. We have examined the
robustness of positive spanning under vector erasures
and its relationship with exactness and cardinality
restrictions, together with approximation properties
of positive Parseval frames.

Several questions remain open for further
investigation. One natural direction is to extend
the present framework to infinite-dimensional
Hilbert spaces, where positive spanning and the
existence of strictly positive relations require
additional care. Another interesting problem is
to investigate positive dual frames under further
constraints, such as prescribed norms, sparsity,
or optimal reconstruction error in the presence
of erasures. It would also be worthwhile to
explore positive frames generated by more general
group actions and to understand how the geometry
of their relation spaces interacts with the
underlying representation. We hope that the
results developed here provide a useful starting
point for these investigations and further
connections between frame theory, convex geometry,
and constrained reconstruction.

\section{Acknowledgment}

\noindent
The first author is grateful to the Mohapatra Family Foundation and the
College of Graduate Studies at the University of Central Florida for
supporting his postdoctoral fellowship during the course of this research.
The second author acknowledges partial support from a Simons Travel Award.

\section{Declaration} 
\noindent Ethical approval: Not applicable.\\
Data Availability Statement: Not Applicable.\\
Conflict of Interest: Not Applicable.\\
Clinical Trial: Not applicable.


\end{document}